  % FONTS

  \newcount\fontset
  \fontset=1
  \def\dualfont#1#2#3{\font#1=\ifnum\fontset=1 #2\else#3\fi}

  \dualfont\bbfive{bbm5}{cmbx5}
  \dualfont\bbseven{bbm7}{cmbx7}
  \dualfont\bbten{bbm10}{cmbx10}
  \font \eightbf = cmbx8
  \font \eighti = cmmi8 \skewchar \eighti = '177
  \font \eightit = cmti8
  \font \eightrm = cmr8
  \font \eightsl = cmsl8
  \font \eightsy = cmsy8 \skewchar \eightsy = '60
  \font \eighttt = cmtt8 \hyphenchar\eighttt = -1
  \font \msbm = msbm10
  \font \sixbf = cmbx6
  \font \sixi = cmmi6 \skewchar \sixi = '177
  \font \sixrm = cmr6
  \font \sixsy = cmsy6 \skewchar \sixsy = '60
  \font \tensc = cmcsc10
  
  \font \titlefont = cmr7 scaled 2000
  \scriptfont \bffam = \bbseven
  \scriptscriptfont \bffam = \bbfive
  \textfont \bffam = \bbten

  \font\rs=rsfs10 % For Lin(H)

  \newskip \ttglue

  \def \eightpoint {\def \rm {\fam0 \eightrm }%
  \textfont0 = \eightrm
  \scriptfont0 = \sixrm \scriptscriptfont0 = \fiverm
  \textfont1 = \eighti
  \scriptfont1 = \sixi \scriptscriptfont1 = \fivei
  \textfont2 = \eightsy
  \scriptfont2 = \sixsy \scriptscriptfont2 = \fivesy
  \textfont3 = \tenex
  \scriptfont3 = \tenex \scriptscriptfont3 = \tenex
  \def \it {\fam \itfam \eightit }%
  \textfont \itfam = \eightit
  \def \sl {\fam \slfam \eightsl }%
  \textfont \slfam = \eightsl
  \def \bf {\fam \bffam \eightbf }%
  \textfont \bffam = \eightbf
  \scriptfont \bffam = \sixbf
  \scriptscriptfont \bffam = \fivebf
  \def \tt {\fam \ttfam \eighttt }%
  \textfont \ttfam = \eighttt
  \tt \ttglue = .5em plus.25em minus.15em
  \normalbaselineskip = 9pt
  \def \MF {{\manual opqr}\-{\manual stuq}}%
  \let \sc = \sixrm
  \let \big = \eightbig
  \setbox \strutbox = \hbox {\vrule height7pt depth2pt width0pt}%
  \normalbaselines \rm }

  % HEADER

  \def \Headlines #1#2{\nopagenumbers
    \voffset = 2\baselineskip
    \advance \vsize by -\voffset
    \headline {\ifnum \pageno = 1 \hfil
    \else \ifodd \pageno \tensc \hfil \lcase {#1} \hfil \folio
    \else \tensc \folio \hfil \lcase {#2} \hfil
    \fi \fi }}

  \def \Title #1{\vbox{\baselineskip 20pt \titlefont \noindent #1}}

  \def \Date #1 {\footnote {}{\eightit Date: #1.}}

  \def \Authors #1{\bigskip \bigskip \noindent #1}

  \long \def \Addresses #1{\begingroup \eightpoint \parindent0pt
\medskip #1\par \par \endgroup }

  \long \def \Abstract #1{\begingroup \eightpoint
  \bigskip \bigskip \noindent
  {\sc ABSTRACT.} #1\par \par \endgroup }

  % CONTROL SEQUENCES

  \def \lcase #1{\edef \auxvar {\lowercase {#1}}\auxvar }

  \def \goodbreak {\vskip0pt plus.1\vsize \penalty -250 \vskip0pt
plus-.1\vsize }

  \newcount \secno \secno = 0
  \newcount \stno

  \def \seqnumbering {\global \advance \stno by 1
    \number \secno .\number \stno }

  \def \label #1{\def\localvariable {\number \secno
    \ifnum \number \stno = 0\else .\number \stno \fi }\global \edef
    #1{\localvariable }}

  \def\section #1{\global\def\SectionName{#1}\stno = 0 \global
\advance \secno by 1 \bigskip \bigskip \goodbreak \noindent {\bf
\number \secno .\enspace #1.}\medskip \noindent \ignorespaces}

  \long \def \sysstate #1#2#3{\medbreak \noindent {\bf \seqnumbering
.\enspace #1.\enspace }{#2#3\vskip 0pt}\medbreak }
  \def \state #1 #2\par {\sysstate {#1}{\sl }{#2}}
  \def \definition #1\par {\sysstate {Definition}{\rm }{#1}}

  \def \proof {\medbreak \noindent {\it Proof.\enspace }}
  \def \proofend {\ifmmode \eqno \square \else \hfill \square
\looseness = -1 \medbreak \fi }

  \def\=#1{\buildrel (#1) \over =}

  \def\Item #1{\smallskip \item {#1}}
  \newcount \zitemno \zitemno = 0
  \def\izitem {\zitemno = 0}
  \def\zitem {\global \advance \zitemno by 1 \Item {{\rm(\romannumeral
\zitemno)}}}

  \newcount \footno \footno = 1
  \newcount \halffootno \footno = 1
  \def\footcntr {\global \advance \footno by 1
  \halffootno =\footno
  \divide \halffootno by 2
  $^{\number\halffootno}$}
  \def\fn#1{\footnote{\footcntr}{\eightpoint#1}}

  % STANDARD DEFINITIONS

  \def \({\left (}
  \def \){\right )}
  \def \[{\left \Vert }
  \def \]{\right \Vert }
  \def \*{\otimes }
  \def \+{\oplus }
  \def \:{\colon }
  \def \<{\left \langle }
  \def \>{\right \rangle }
  \def \text #1{\hbox {\rm #1}}
  \def \and {\hbox {,\quad and \quad }}
  
  \def \calcat #1{\,{\vrule height8pt depth4pt}_{\,#1}}

  \def \crossproduct {{\hbox {\msbm o}}}
  
  \def \for #1{,\quad \forall\,#1}
  \def \inv {^{-1}}
  
  \def \square {\hbox {$\sqcap \!\!\!\!\sqcup $}}
  \def \stress #1{{\it #1}\/}
  
  \def \x {\times }
  \def \|{\Vert }
  \def \inv {^{-1}}

  % REFERENCE CONTROL

  \def\cite #1{{\rm [\bf #1\rm ]}}
  \def\scite #1#2{\cite{#1{\rm \hskip 0.7pt:\hskip 2pt #2}}}
  \def\lcite #1{(#1)}
  \def\fcite#1{#1}
  \def\bibitem#1#2#3#4{\smallskip \item {[#1]} #2, ``#3'', #4.}

  \def \references {
    \begingroup
    \bigskip \bigskip \goodbreak
    \eightpoint
    \centerline {\tensc References}
    \nobreak \medskip \frenchspacing }

  \def\tsigma{{\hat \sigma}}
  \def\tS{{\hat S}}
  \def\tK{{\hat \K}}
  \def\te{{\hat e}}
  \def\tgamma{{\hat \gamma}}
  \def\tG{{\hat G}}
  \def\tU{{\moveaccent{\hat}{{\curly U}}{2}}}
  \def\tP{{\hat P}}

  \def\csigma{\sigma}
  \def\cS{S}
  \def\cK{\K}
  \def\ce{e}
  \def\cgamma{\gamma}
  \def\cG{G}
  \def\cU{{\curly U}}
  \def\cMCP{\cU\mathop{\crossproduct_{\b}}\N}

  \def\mS{{\check S}}
  \def\mK{{\check\K}}
  \def\me{{\check e}}
  \def\mgamma{{\check \gamma}}
  \def\mF{{\check F}}
  \def\ma{{\check \a}}
  \def\mtr{{\check \Tr}}
  \def\mU{{\moveaccent{\check}{{\curly U}}{2}}}

  \def\moveaccent#1#2#3{ \def \ofs{#3pt} \kern \ofs #1 {\kern -\ofs #2
    \kern \ofs} \kern -\ofs}
  \def \&#1{#1 $$$$ #1}
  \def\pilar#1{\vrule height #1pt width 0pt}
  \def\kms{KMS$_\b$}
  \def\iota{I}
  \def\tensor{A\*_{\R_n}\!\!A}
  \def\comm#1#2{C_{#2}(#1)}
  \def\({\Big(}
  \def\){\Big)}
  \def\ind{{\sl ind}}
  \def\a{\alpha}
  \def\b{\beta}
  \def\numb{N}
  \def\K{{\cal K}}
  \def\R{{\cal R}}
  \def\E{{\cal E}}
  
  \def\Ker{{\sl Ker}}
  \def\N{{\bf N}}
  \def\Z{{\bf Z}}
  \def\Real{{\bf R}}
  \def\C{{\bf C}}
  \def\Tr{{\cal L}}
  \def\T{{\cal T}}
  \def\map#1{\mathrel{\buildrel #1 \over \longrightarrow}}

  \def\curly#1{\hbox{\rs #1\/}}
  \def\CP{A\mathop{\crossproduct_{\a,\Tr}} {\bf N}}
  \def\TCP{{\curly T}(A,\a,\Tr)}
  \def\MCP{\mU\mathop{\crossproduct_{\b}}\N}
  \def\compos{\mathop{\raise 1pt \hbox{$\scriptscriptstyle \circ$}}}
  \def\Lin {{\curly L}}
  \def\Zenter {{\curly Z}}
  \def\nondeg{non-degenerate}
  \font\forsmallbullet=cmsy10 scaled 500
  \def\smallbullet{\raise 1pt \hbox{\forsmallbullet }}

  % REFERENCES

  \def\Endo{E2}
  \def\Takesaki{T}
  \def\newpim{E1}
  \def\Ped{P}
  \def\Murphy{M}
  \def\Stacey{S}
  \def\CuntzTwo{C}
  \def\Watatani{W}
  \def\Kosaki{K}
  \def\Jones{J}
  \def\RueleOne{R1}
  \def\RueleTwo{R2}

  \def\titletext{Crossed-Products by Finite Index Endomorphisms and KMS
states}

  \Headlines {\titletext} {Ruy Exel}

  \Title {\titletext}\footnote{\null}
  {\eightrm 2000 \eightsl Mathematics Subject Classification:
  \eightrm 
  46L55.% Noncommutative dynamical systems
  }

  \Date {May 23, 2001}

  \Authors
  {Ruy Exel\footnote{*}{\eightrm Partially supported by CNPq.}}

  \Addresses
  {Departamento de Matem\'atica\par
  Universidade Federal de Santa Catarina\par
  88040-900 Florian\'opolis SC\par
  BRAZIL\vskip 4pt
  E-mail: exel@mtm.ufsc.br}

  \Abstract {Given a unital C*-algebra $A$, an injective endomorphism
$\a\colon A \to A$ preserving the unit, and a conditional expectation
$E$ from $A$ to the range of $\a$ we consider the crossed-product of
$A$ by $\a$ relative to the transfer operator $\Tr=\a\inv E$.  When
$E$ is of index-finite type we show that there exists a conditional
expectation $G$ from the crossed-product to $A$ which is unique under
certain hypothesis.
  We define a ``gauge action'' on the crossed-product algebra in terms
of a central positive element $h$ and study its KMS states.  The main
result is: if $h>1$ and $E(ab)=E(ba)$ for all $a,b\in A$ (e.g.~when
$A$ is commutative) then the {\kms} states are precisely those of the
form $\psi = \phi\compos G$, where $\phi$ is a trace on $A$ satisfying
the identity
  $$
  \phi(a) = \phi(\Tr(h^{-\b}\ind(E)a)),
  $$
  where
  $\ind(E)$ is the Jones-Kosaki-Watatani index  of $E$.
  }

  \section{Introduction}
  In \cite{\Endo} we have introduced the notion of the crossed-product
of a C*-algebra $A$ by a *-endomorphism $\a$, a construction which
also depends on the choice of a \stress{transfer operator}, that is a
positive continuous linear map $\Tr:A\to A$ such that
$\Tr\big(\a(a)b\big) = a\Tr(b)$, for all $a,b\in A$.  In the present
work we treat the case in which $\a$ is a monomorphism (injective
endomorphism) and $\Tr$ is given by $\Tr=\a\inv\compos E$, where $E$
is a conditional expectation onto the range of $\a$.

  The first of our main results (Theorem \fcite{4.12}) is the solution to
a problem posed in \cite{\Endo}: we prove that the canonical mapping
of $A$ into $\CP$ is injective.  The main technique used to accomplish
this is based on the celebrated ``Jones basic construction''
\cite{\Jones}, as adapted to the context of C*-algebras by Watatani
\cite{\Watatani}.  In order to briefly describe this technique
consider, for each $n\in\N$, the conditional expectation onto the
range of $\a^n$ given by
  $$
  \E_n = \a^n \underbrace{(E \a\inv) \ldots (E \a\inv)}_{n\rm\;times}E.
  $$
  Let $\cK_n$ be the \stress{C*-basic construction}
\scite{\Watatani}{Definition 2.1.10} associated to $\E_n$ and let
$\ce_n$ be the standard projection as in \scite{\Watatani}{Section
2.1}.  We are then able to find a simultaneous representation of all
of the $\cK_n$ in a fixed C*-algebra.  Since $\K_0=A$ we then have
that $A$ is also represented there and we find that $\ce_{n+1}\leq
\ce_n$ for all $n$.

Letting $\cU$ be the C*-algebra generated by the union of all the
$\cK_n$ we construct an endomorphism $\b$ of $\cU$ which is not quite
an extension of $\a$ but which satisfies $\b(a) = \a(a)\ce_1$ and
$\b(\ce_n) = \ce_{n+1}$.

It turns out that the range of $\b$ is a hereditary subalgebra of
$\cU$ and hence we may form the crossed-product $\cMCP$.  We then
prove that $\CP$ is isomorphic to $\cMCP$ (Theorem \fcite{6.5}).  

Crossed-products by endomorphisms with hereditary range are much
easier to understand.   In particular it is known that $\cU$  embeds
injectively in $\cMCP$ and hence we deduce that $A$ is faithfully
represented in $\CP$ as already mentioned.

Another consequence of the existence of an isomorphism between $\cMCP$
and $\CP$ is that we get a rather concrete description of the
structure of $\CP$ and, in particular, of the fixed point subalgebra
for the (scalar) gauge action, namely the action of the circle on
$\CP$ given by
  $$
  \cgamma_z(\cS ) = z\cS
  \and
  \cgamma_z(a) = a
  \for a\in A
  \for z\in S ^1,
  $$
  where $S$ is the standard isometry in $\CP$.
  That fixed-point algebra, if viewed from the point of view of $\cMCP$, is
well known to be exactly $\cU$ (see \scite{\Murphy}{4.1}).

Given an action of the circle on a C*-algebra there is a standard way
to construct a conditional expectation onto the fixed-point algebra by
averaging the action.  It is therefore easy to construct a conditional
expectation from $\CP$ to $\cU$.  The existence of a conditional
expectation onto $A$, however, is an entirely different matter.

Our second main result (Theorem \fcite{8.9}) is the construction of
such a conditional expectation under the special case in which $E$ is
of index-finite type \scite{\Watatani}{1.2.2}.  Precisely we show that
there is a (unique under certain circumstances) conditional
expectation $\cG : \CP \to A$ such that
  $$
  \cG(a\cS^n\cS^{*m}b) =
  \delta_{nm}a\iota_n\inv b
  \for a,b\in A \for n,m\in\N,
  $$
  where $\delta$ is the Kronecker symbol,
  $$
  \iota_n = 
  \ind(E) \a\big(\ind(E)\big) \ldots \a^{n-1}\big(\ind(E)\big),
  $$
  and $\ind(E)$ is the index of $E$ defined by Watatani
in \scite{\Watatani}{1.2.2}, generalizing earlier work of Jones \cite{\Jones}
and Kosaki \cite{\Kosaki}.

Our third main result (Theorem \fcite{9.6}) is related to the KMS
states on $\CP$ for the one-parameter automorphism group $\csigma$ of
$\CP$ specified by
  $$
  \csigma_t(\cS )=h^{it}\cS  \and
  \csigma_t(a)=a
  \for a\in A,
  $$  
  where $h$ is any self-adjoint element in the center of $A$ such that
$h \geq cI$ for some real number $c>1$.  Under the hypothesis that $E$
is of index-finite type, and hence in the presence of the conditional
expectation $\cG$ above, and also assuming that
$E(ab)=E(ba)$ for all $a,b\in A$ (e.g.~when $A$ is commutative),
we show that all KMS states on $\CP$ factor
through $G$ and are exactly the states $\psi$ on $\CP$ given by $\psi=\phi\compos G$
where $\phi$ is a trace on $A$ such that
  $$
  \phi(a) = \phi\big(\Tr(h^{-\b}\ind(E)a)\big)
  $$
  for all $a\in A$.
  We also show that there are no ground states on $\CP$.

  We conclude with a brief discussion of the case in which $A$ is
commutative and show that the KMS states on $\CP$ are related to
Ruelle's work on Statistical Mechanics \cite{\RueleOne},
\cite{\RueleTwo}.

  A word about our notation: most of the time we will be working
simultaneously with three closely related algebras, namely the
``Toeplitz extension'' $\TCP$, the crossed-product $\CP$, and a
concretely realized algebra $\MCP$.  The features of each of these
algebras will most of the time be presented side by side, e.g.~each
one will contain a distinguished isometry.  In order to try to keep
our notation simple but easy to understand we have chosen to decorate
the notation relative to the first algebra with a ``hat'', the one for
the third with a ``check'', and no decoration at all for $\CP$ which
is, after all, the algebra that we are most interested in.  For
example, the three isometries considered will be denoted $\tS$, $\cS$,
and $\mS$.

  We would like to acknowledge helpful conversations with Marcelo
Viana from which some of the intuition for the present work
developed.

  \section{Crossed products}
  Throughout this section, and most of this work, we will let $A$ be a
unital C*-algebra and
  $
  \a:A\to A
  $
  be an injective *-endomorphism such that $\a(1)=1$.  It is
conceivable that some of our results survive without the hypothesis
that $\a$ be injective but for the sake of simplicity we will stick
to the injective case here.

The range of $\a$, which will play a predominant role in what follows,
will be denoted by $\R$ and we will assume the existence of a
  {\nondeg}\fn{A conditional expectation $E$ is said to be {\nondeg}
when $E(a^*a)=0$ implies that $a=0$.}
  conditional expectation
  $$
  E: A \to \R
  $$
  which will be fixed throughout. As in \scite{\Endo}{2.6} it follows
that the composition
  $\Tr:= \a\inv E$
  is a \stress{transfer operator} in the sense of \scite{\Endo}{2.1},
meaning a positive linear map $\Tr:A\to A$ such that $\Tr\big(\a(a)b\big) =
a\Tr(b)$, for all $a,b\in A$.

According to Definition 3.1 in \cite{\Endo} the ``Toeplitz extension''
$\TCP$ is the universal unital C*-algebra generated by a copy of
$A$ and an element $\tS$ subject to the relations:
  \izitem
  \zitem $\tS a = \a(a)\tS $, and
  \zitem $\tS^*a\tS  = \Tr(a)$,
  \medskip\noindent
  for every $a\in A$.
  As proved in \scite{\Endo}{3.5} the canonical map from $A$ to $\TCP$
is injective so we may and will view $A$ as a subalgebra of $\TCP$.

Observe that, as a consequence of the fact that $\a$ preserves the
unit, we have that $1\in \R$ and hence that
  $
  \Tr(1) = \a\inv(E(1)) = 1.
  $
  It follows that
  $$
  \tS^*\tS  =   \tS^*1\tS  = \Tr(1) = 1,
  $$
  and hence we see that $\tS $ is an isometry.

Following \scite{\Endo}{3.6} a \stress{redundancy} is a pair $(a,k)$
of elements in $\TCP$ such that $k$ is in the closure of $A\tS
\tS^*A$, $a$ is in $A$, and
  $$
  ab\tS  = kb\tS  \for b\in A.
  $$

  \definition
  \scite{\Endo}{3.7}
  The crossed-product of $A$ by $\a$ relative to $\Tr$, denoted by
$\CP$, is defined to be the quotient of $\TCP$ by the closed two-sided
ideal generated by the set of differences $a-k$, for
  all\fn{We should remark that in Definition 3.7 of \cite{\Endo} one
uses only the redundancies $(a,k)$ such that $a\in\overline{A\R A}$.
But, under the present hypothesis that $\a$ preserves the unit, we
have that $1\in \R$ and hence $\overline{A\R A}=A$.}  redundancies
$(a,k)$.
  We will denote by $q$ the canonical quotient map
  $$
  q: \TCP \to \CP,
  %  \eqno {(\seqnumbering)}
  % \label \NotationQuotient
  $$
  and by $\cS$ the image of $\tS$ under $q$.

  \state Lemma
  \label \FormulasForLinearSpan
  Given $n,m,j,k\in\N$ and $a,b,c,d\in A$ let $x,y\in\TCP$ be given by
$x=a\tS^n\tS^{*m}b$ and $y=c\tS^j\tS^{*k}d$.  Then
  $$
  xy = \cases{   a\a^n(\Tr^m(bc))\tS^{n-m+j}\tS^{*k}d, & if $m\leq j$, \cr
  a\tS^n\tS^{*(m-j+k)}\a^k(\Tr^j(bc))d, & if $m \geq j$.}
  $$

  \proof
  If $m\leq j$ one has
  $$
  xy =
  a\tS^n(\tS^{*m}bc\tS^m)\tS^{j-m}\tS^{*k}d =
  a\tS^n\Tr^m(bc)\tS^{j-m}\tS^{*k}d =
  a\a^n(\Tr^m(bc))\tS^{n+j-m}\tS^{*k}d.
  $$
  On the other hand, if $m \geq j$ one has
  $$
  xy =
  a\tS^n\tS^{*(m-j)}(\tS^{*j}bc\tS^j)\tS^{*k}d =
  a\tS^n\tS^{*(m-j)}\Tr^j(bc)\tS^{*k}d =
  a\tS^n\tS^{*(m-j+k)}\a^k(\Tr^j(bc))d.
  \proofend
  $$

  As a consequence we have:

  \state Proposition
  \label \LinearSpan
  $\TCP$ is the closed linear span of the set $X= \{a\tS^n\tS^{*m}b: a,b\in A,\
n,m\in \N\}$.

  \proof
  By \lcite{\FormulasForLinearSpan} we see that the linear span of $X$
is an algebra.  Since it is also self-adjoint and contains $A \cup
\{\tS \}$ the result follows.
  \proofend

  There are may results for $\TCP$ which yield similar results for
$\CP$ simply by passage to the quotient, such as
\lcite{\FormulasForLinearSpan} and \lcite{\LinearSpan}.  Most often we
will not bother to point these out unless it is relevant to our
purposes that we do so.

  \section{Gauge action}
  In this section we will describe certain one-parameter automorphism groups
of $\TCP$ and $\CP$ relative to which we will later study KMS states.

  \state Proposition
  \label \AutomorphismFromU
  Given a unitary element $u$ in $\Zenter(A)$ (the center of $A$) there exists a
unique automorphism $\tsigma_u$ of $\TCP$ such that
  $$
  \tsigma_u(\tS )=u\tS  \and
  \tsigma_u(a)=a
  \for a\in A.
  $$
  Moreover
$\tsigma_u(\Ker(q))=\Ker(q)$ and hence $\tsigma_u$ drops to the quotient
providing an automorphism $\csigma_u$ of $\CP$ such that
  $$
  \csigma_u(\cS )=q(u)\cS  \and
  \csigma_u(a)=a
  \for a\in A.
  $$
  If $v$ is another unitary element in $\Zenter(A)$ then $\tsigma_u
\tsigma_v=\tsigma_{uv}$.

  \proof
  Let $S_u = u\tS $ and observe that for every $a$ in $A$ one has
  $$
  S_u a = u\tS a = u\a(a)\tS  = \a(a)u\tS  = \a(a)S_u,
  $$
  and
  $$
  S_u^*aS_u =   \tS^*u^*au\tS  =  \tS^*a\tS  = \Tr(a).
  $$
  From the universal property of $\TCP$ it follows that there exists a
unique *-homomorphism $\tsigma_u:\TCP\to\TCP$ such that
$\tsigma_u(a)=a$, for all $a\in A$, and $\tsigma_u(\tS )=S_u$.  Given
$v$ as above notice that
  $$
  \tsigma_u \tsigma_v(\tS )=
  \tsigma_u(v\tS ) =
  \tsigma_u(v) \tsigma_u(\tS ) =
  vu\tS  = uv\tS  = \tsigma_{uv}(\tS ),
  $$
  and that $\tsigma_u \tsigma_v(a)=a$, for all $a\in A$.  Thus
$\tsigma_u \tsigma_v=\tsigma_{uv}$.  It follows that $\tsigma_{u\inv}$ is
the inverse of $\tsigma_u$ and hence $\tsigma_u$ is an automorphism.

Let $(a,k)$ be a redundancy.  Then
  $$
  \tsigma_u(k) \in \tsigma_u(\overline{A\tS \tS^*A}) =
  \overline{Au\tS \tS^*u^*A} =
  \overline{A\tS \tS^*A}.
  $$
  For every $b$ in $A$ we have
  $$
  \tsigma_u(k)b\tS  =
  \tsigma_u(k)bu\inv u\tS  =
  \tsigma_u(kbu\inv \tS ) =
  \tsigma_u(abu\inv \tS ) =
  ab\tS ,
  $$
  so $(a,\tsigma_u(k))$ is also a redundancy and it follows that
$\tsigma_u(a-k)\in\Ker(q)$ and hence that
  $\tsigma_u(\Ker(q))\subseteq \Ker(q)$.  Since the same holds for
$\tsigma_{u\inv}$ we have that
  $\Ker(q) \subseteq \tsigma_u(\Ker(q))$.
  \proofend

  Let $h\in\Zenter(A)$ be a self-adjoint element such that $h \geq cI$
for some real number $c>0$.  For every $t\in\Real$ we have that
$h^{it}$ is a unitary in $\Zenter(A)$ and hence defines an
automorphism $\tsigma_{h^{it}}$ by \lcite{\AutomorphismFromU} which we
will denote by $\tsigma_t^h$.  Again by \lcite{\AutomorphismFromU} we
have that $\tsigma^h_t \tsigma^h_s = \tsigma^h_{t+s}$ so that
$\tsigma^h$ is a one-parameter automorphism group of $\TCP$ which is
clearly strongly continuous.

  \definition
  \label \DefineGaugeAction
  Both the action $\tsigma^h$ defined above and the action $\csigma^h$
of $\Real$ on $\CP$ obtained by passing $\tsigma^h$ to the quotient
will be called the \stress{gauge action associated to $h$}.

When $h$ is taken to be Neper's number $e$ we have that
$\tsigma_t^h(\tS ) = e^{it}\tS$, so the gauge action is periodic with
period $2\pi$ and hence defines an action $\tgamma$ of the unit circle
on $\TCP$ such that
  $$
  \tgamma_z(\tS ) = z\tS
  \and
  \tgamma_z(a) = a
  \for a\in A
  \for z\in S ^1.
  $$

  \definition
  \label \ScalarGaugeAction
  Both the action $\tgamma$ defined above and the action $\cgamma$ of
the circle group on $\CP$ obtained by passing $\tgamma$ to the
quotient will be called the \stress{scalar gauge action}.

We will later be interested in the fixed point algebra for the scalar
gauge action so the following result will be useful:

  \state Proposition
  \label \PreFixedPointAlgebras
  Let $B$ be a C*-algebra with a strongly continuous action $\gamma$
of the circle group.  Suppose that $B$ is the closed linear span of a
set $\{x_i: i\in I\}$ such that for every $i\in I$ there exists
$n_i\in\Z$ such that 
  $
  \gamma_z(x_i) = z^{n_i} x_i
  $
  for all $z\in\C$.  Then the fixed point algebra for $\gamma$ is the
closed linear span of $\{x_i: n_i=0\}$.

  \proof
  It is well known that the map $P:B\to B$ given by 
  $$
  P(a)= \int_{S^1} \gamma_z (a)\, dz 
  $$  
  is a conditional expectation onto the fixed point algebra for
$\gamma$.  By direct computation it is easy to see that $P(x_i) = 0$
when $n_i\neq 0$ and  $P(x_i) = x_i$
when $n_i= 0$.

Given a fixed point $b$ and $\varepsilon>0$ let
$\{\lambda_i\}_i$ be a family of scalars with finitely many nonzero
elements such that
  $
  \| b - \sum_{i\in I} \lambda_i x_i \| < \varepsilon.
  $
  It follows that 
  $$
  \Big\Vert b - \sum_{
  {\buildrel {\scriptstyle i\in I} \over {n_i=0}}
  } \lambda_i x_i \| =
  \Big\|P\Big(b - \sum_{i\in I} \lambda_i x_i\Big) \Big\| \leq
  \Big\|b - \sum_{i\in I} \lambda_i x_i \Big\| < 
  \varepsilon.
  $$
  Therefore $b\in\overline{\sl span}\{x_i: n_i=0\}$.
  \proofend

  \state  Corollary 
  \label \FixedPointAlgebras
  The fixed point subalgebra of $\TCP$ (resp.~$\CP$)
for the scalar gauge action $\tgamma$
(resp.~$\cgamma$) is the closed linear span of the set of elements
$a\tS^n\tS^{*n}b$ (resp.~$a\cS^n\cS^{*n}b$) for all $a,b\in A$ and
$n\in\N$.

  \section{Conditional Expectations and Hilbert Modules}
  In this section we will describe certain conditional expectations
and certain Hilbert modules which will be used in later sections.
  For every $n\in\N$ we shall let $\R_n$ denote the range of $\a^n$.
Therefore $\R_0=A$, $\R_1=\R$, and the $\R_n$ form a descending chain
of closed *-subalgebras of $A$
  $$
  A = \R_0 \supseteq \R_1 \supseteq \R_2 \supseteq \ldots.
  $$
  Clearly  each $\R_n$ is isomorphic to $A$ under $\a^n$.
For each $n\in \N$ consider the map
  $$
  E_n : \R_n \to \R_{n+1}
  $$
  given by $E_n = \a^n E \a^{-n}$.  It is elementary to verify that
each $E_n$ is a {\nondeg} conditional expectation.
  % and that $E_{n+1} = \a E_n\a\inv$ for all $n$.
  Likewise, for each $n\in\N$, the composition
  $$
  A \map{E_0} \R_1 \map{E_1} \R_2 \map{E_2} \ \cdots\ \map{E_{n-1}} \R_n
  $$
  is a {\nondeg} conditional expectation onto $\R_n$, which we denote
by $\E_n$.  By default we let $\E_0$ be the identity map on $A$ and it
is clear that $\E_1 = E_0$ = $E$.

For future use it is convenient to record the following elementary
facts:

  \state Proposition
  \label \FormulasForFn
  For every $n\in \N$ one has that
  \izitem
  \zitem $\E_{n+1} = E_n \E_n = \a \E_n \a\inv E$,
  \zitem $\E_{n+1}\E_n = \E_n  \E_{n+1} = \E_{n+1}$.

  We now need to use a simple construction from the theory of Hilbert
modules:
  let $B$ be any C*-algebra and let $C\subseteq B$ be a
sub-C*-algebra.  Also let
  $
  E: B \to C
  $
  be a {\nondeg} conditional expectation.  Given a right Hilbert
$B$--module $M$ (with inner--product $\<\cdot\thinspace,\cdot\>$)
one gets a $C$--valued inner--product on $M$ defining
  $$
  \<x,y\>_C = E\big(\<x,y\>\big) \for x,y\in M.
  $$
  We shall denote the Hilbert $C$--module obtained by
completing $M$ under the norm $\|x\|_C = \|\<x,y\>_C\|^{1/2}$ by $M_C$.

We plan to apply this construction in order to obtain a sequence
$\{M_n\}_{n\in \N}$, where each $M_n$ is a Hilbert $\R_n$-module as
follows: let $M_0=A$ viewed as a right Hilbert $A$--module under the
obvious right module structure and inner--product given by $\<a,b\> =
a^*b$, for all $a$ and $b$ in $A$.

  Once $M_n$ is constructed let $M_{n+1}$ be the $\R_{n+1}$--module
obtained by applying the procedure described above to $M_n$ and the conditional
expectation $E_n$.  For simplicity we let
$\<\cdot\thinspace,\cdot\>_n$ denote the inner--product on $M_n$ and
by $\|\cdot\|_n$ the associated norm.  By construction we have that
  $$
  A = M_0 \subseteq M_1 \subseteq M_2 \subseteq
\cdots
  $$
  where the inclusion maps are continuous and each $M_n$ is a dense
subset of $M_{n+1}$ (with respect to $\|\cdot\|_{n+1}$).  It follows
that $A$ is dense in each $M_n$ and it is convenient to observe that
  $$
  \<a,b\>_n = E_{n-1}\cdots E_1 E_0(a^*b) =  \E_n(a^*b) \for a,b\in A.
  $$

  \state Proposition
  For every $n\in\N$
  there exists an isometric complex-linear map
  $\ma_n : M_n \to M_{n+1}$
  such that
  $\ma_n(a) = \a(a)$  for all
  $a\in A$.

  \proof
  Given $a\in A$ we have
  $$
  \<\a(a),\a(a)\>_{n+1} =
  \E_{n+1}(\a(a^*a)) =
  \a \E_n \a\inv E(\a(a^*a)) =
  \a \E_n (a^*a) =
  \a\big(\<a,a\>_n\big).
  $$
  This implies that  $\|\a(a)\|_{n+1} = \|a\|_n$ from where the
conclusion easily follows.
  \proofend

  \state Proposition
  For every $n\in\N$
  there exists a contractive complex-linear map
  $\mtr_n : M_{n+1} \to M_n$
  such that
  $\mtr_n(a) = \Tr(a)$ for all
$a\in A$.

  \proof
  Using the well known fact that $E(a^*)E(a) \leq E(a^*a)$  (plug
$x:=a-E(a)$ in ``$E(x^*x)\geq0$'' in order to prove it) we have that
  $$
  \<\Tr(a),\Tr(a)\>_n =
  \E_n\big(\a\inv(E(a)^*E(a))\big) \leq
  \E_n \a\inv E(a^*a) \&=
  \a\inv\a \E_n\a\inv E(a^*a) =
  \a\inv \E_{n+1}(a^*a) =
  \a\inv\big(\<a,a\>_{n+1}\big).
  $$
  We then have that
  $
  \|\Tr(a)\|_n\leq \|a\|_{n+1}
  $
  from where one easily deduces the existence of $\mtr_n$.
  \proofend

  From now on we will denote by $\Lin(M_n)$ the C*-algebra of all
adjointable operators on $M_n$.

  \state Proposition
  \label \IntroductionOfEn
  For every $n\in\N$
  there exists a self-adjoint idempotent $\me_n\in\Lin(M_n)$
  such that
  $\me_n(a) = \E_n(a)$ for all $a\in A$.

  \proof
  For $a\in A$ we have
  $$
  \<\E_n(a),\E_n(a)\>_n =
  \E_n\big(\E_n(a^*)\E_n(a)\big) =
  \E_n(a^*)\E_n(a) \leq
  \E_n(a^*a) =
  \<a,a\>_n.
  $$
  Therefore $\|\E_n(a)\|_n\leq\|a\|_n$ and hence
  the correspondence $a \mapsto \E_n(a)$ extends to a bounded
linear map
  $
  \me_n : M_n \to M_n
  $.
  For $a,b\in A$ we have that
  $$
  \<\me_n(a),b\>_n = \E_n(\E_n(a)^*b) = \E_n(a^*)\E_n(b) =
  \E_n(a^*\E_n(b)) =
  \<a,\me_n(b)\>_n.
  $$
  By continuity it follows that $\<\me_n(\xi),\eta\> =
\<\xi,\me_n(\eta)\>$ for all $\xi,\eta\in M_n$ so that $\me_n$ is in fact
an adjointable operator on $M_n$.  The remaining assertions are now
easy to prove.
  \proofend

It should be remarked that $\me_n$ is precisely the projection
introduced in \scite{\Watatani}{Section 2.1} relative to the
conditional expectation $\E_n$.  Therefore $\overline{A\me_n A}$ is
the associated \stress{reduced C*-basic construction}
\scite{\Watatani}{Definition 2.1.2}.  We will soon have more to say
about this.

  \state Proposition
  \label \VariousFormulas
  For every $n\in\N$
  \izitem
  \zitem $\mtr_n \ma_n$ is the identity on $M_n$,
  \zitem $\<\ma_n(\xi),\eta\>_{n+1} =
\a\big(\<\xi,\mtr_n(\eta)\>_n\big)$, for all $\xi\in M_n$ and $\eta
\in M_{n+1}$, and
  \zitem $\ma_n \me_n \mtr_n = \me_{n+1}$.

  \proof
  With respect to (i) we have for all $a\in A$ that
  $$
  \mtr_n \ma_n (a) =
  \Tr\a(a) =
  \a\inv E \a(a) =
  %  \a\inv \a(a) =
  a,
  $$
  so the conclusion follows by continuity.
  As for (ii) one has for all $a,b\in A$ that
  $$
  \<\a(a),b\>_{n+1}  =
  \E_{n+1}(\a(a^*)b) =
  \a \E_n \a\inv E(\a(a^*)b) =
  \a \E_n \a\inv( \a(a^*)E(b)) \&=
  \a \E_n( a^*\Tr(b)) =
  \a\big(\<a,\Tr(b)\>_n\big).
  $$
  Speaking of (iii), fix $a\in A$  and notice that
  $$
  \ma_n \me_n \mtr_n(a) =
  \a \E_n \Tr (a) =
  \a \E_n \a\inv E (a) =
  \E_{n+1} (a) =
  \me_{n+1} (a).
  \proofend
  $$

  \state Proposition
  \label \DefineBeta
  For each $n\in\N$ the map
  $$
  \b_n : T \in  \Lin(M_n) \longmapsto \ma_n T \mtr_n \in \Lin(M_{n+1})
  $$
  is a well defined *-monomorphism of C*-algebras.

  \proof For each $T \in \Lin(M_n)$ it is clear that $\ma_n T \mtr_n$
is a bounded complex-linear map on $M_n$.  Given $\xi,\eta\in M_{n+1}$
notice that
  $$
  \<\b_n(T)\xi,\eta\>_{n+1} =
  \<\ma_n T \mtr_n(\xi),\eta\>_{n+1} =
  \a\big(\<T \mtr_n(\xi),\mtr_n(\eta)\>_n\big) \&=
  \a\big(\<\mtr_n(\xi),T^*\mtr_n(\eta)\>_n\big) =
  \<\xi,\ma_n T^*\mtr_n(\eta)\>_{n+1} =
  \<\xi,\b_n(T^*)\eta\>_{n+1}.
  $$
  This proves that $\b_n(T)$ is an adjointable operator on $M_{n+1}$
with $\b_n(T)^*=\b_n(T^*)$.  So $\b_n$ is a well defined linear map
from $\Lin(M_n)$ to $\Lin(M_{n+1})$ which moreover respects the
involution.  Given $T,S\in \Lin(M_n)$ we have that
  $$
  \b_n(T)  \b_n(S) =
  \ma_n T \mtr_n   \ma_n S \mtr_n =
  \ma_n T S \mtr_n =
  \b_n(TS),
  $$
  proving that $\b$ is a *-homomorphism.
  Suppose that  $T\in \Lin(M_n)$  is such that  $\b_n(T)=0$.  Then
  $$
  0 =
  \mtr_n \b_n(T) \ma_n =
  \mtr_n  \ma_n T \mtr_n \ma_n =
  T.
  $$
  Therefore  $\b_n$ is injective.
  \proofend

  We now need another result from the theory of Hilbert modules.

  \state Lemma
  Under the assumption that $E:B\to C$ is a {\nondeg} conditional
expectation, and $M$ is a right Hilbert $B$--module,
  there exists an injective *-homomorphism
  $$
  \Phi: \Lin_B(M) \to \Lin_C(M_C),
  $$
  such that $\Phi(T)(\xi) = T(\xi)$, for all $T\in \Lin_B(M)$, and all
$\xi\in M$.

  \proof
  Let $T\in \Lin_B(M)$.  Since $T^*T \leq \|T\|^2$ one has for
all $\xi\in M$ that
  $$
  \<T(\xi),T(\xi)\> =
  \<T^*T(\xi),\xi\> \leq
  \|T\|^2\<\xi,\xi\>.
  $$
  Applying $E$ to the above inequality yields
  $$
  \<T(\xi),T(\xi)\>_C \leq
  \|T\|^2\<\xi,\xi\>_C,
  $$
  and hence we conclude that $\|T(\xi)\|_C \leq \|T\|\,\|\xi\|_C$, so
that $T$ is bounded with respect to $\|\cdot\|_C$ and hence extends to
a bounded linear map $\Phi(T)$ on $M_C$.  We leave it for the reader to
verify that $\Phi(T)$ indeed belongs to $\Lin_C(M_C)$ and that the
correspondence $T \mapsto \Phi(T)$ is a *-homomorphism.

Given that $\Phi(T)$ is an extension of $T$ it is clear that
$\Phi(T)\neq 0$ when $T\neq0$, so  that $\Phi$ is injective.
  \proofend

Applying the above result to the present situation
we obtain an inductive sequence  of C*-algebras
  $$
  A = \Lin(M_0) \map {\Phi_0} \Lin(M_1) \map {\Phi_1} \cdots \map
{\Phi_{n-1}}  \Lin(M_n) \map {\Phi_n} \cdots
  $$
  We temporarily denote the inductive limit of this sequence (in the
category of C*-algebras) by $B$.

  \definition   We will denote by $\mU$ the
sub-C*-algebra of $B$ generated by $A \cup \{\me_n\: n\in\N\}$.

The following provides some useful information on the algebraic
structure of $\mU$ (see also \cite{\Jones}, \cite{\Watatani}):

  \state Proposition
  \label \AlgebraicStructureofF
  For every $n\in \N$ one has
  \izitem
  \zitem $\me_{n+1} \leq \me_n$,
  \zitem $\me_n a \me_n = \E_n(a) \me_n = \me_n \E_n(a)$,
  \zitem the linear span of the set
  $\{a\,\me_n b : n\in \N, \ a,b\in A\}$ is dense in $\mU$.

  \proof In order to prove (i) we need to verify that
  $\me_{n+1}\me_n = \me_{n+1}$ or, more precisely, that
  $\me_{n+1}\Phi_n(\me_n) = \me_{n+1}$ as operators on $M_{n+1}$.
  Given any $a\in A \subseteq
M_{n+1}$ we have
  $$
  \me_{n+1}\Phi_n(\me_n)(a) =
  \me_{n+1}\me_n(a) =
  \E_{n+1}\E_n(a) =
  \E_{n+1}(a) =
  \me_{n+1}(a) .
  $$
  For all $b\in A \subseteq M_n$ we have
  $$
  \me_n a \me_n(b) =
  \E_n (a \E_n(b)) =
  \E_n (a) \E_n(b) =
  \E_n (a) \me_n(b).
  $$
  This proves that 
$\me_n a \me_n = \E_n(a) \me_n$.  Taking adjoints it follows that 
$\me_n a \me_n =  \me_n\E_n(a)$ also. Now let $\mU_0$ be the linear span of the set
described in (iii).  We claim that it is a *-subalgebra of
$\mU$. Clearly $\mU_0$ is self-adjoint so we are left with the task of
checking it to be closed under multiplication.  In order to see it let
$a,b,c,d\in A$, and $n,m\in \N$.   We then have
  $$
  (a\me_n b)(c \me_m d) =
  a( \me_n b c \me_n) \me_m d = \cdots
  $$
  where we are assuming, without loss of generality, that $m\geq n$
and hence that $\me_m \leq \me_n$ by (i).  Using  (ii) we conclude that the
above equals
  $$
  \cdots =
  a\E_n(b c) \me_n \me_m d =
  a\E_n(b c) \me_m d,
  $$
  which is seen to belong to $\mU_0$.  This proves our claim.
  It is evident that $A\subseteq \mU_0$ (because $\me_0=1$), and that
$\me_n\in \mU_0$ for all $n$.   Since $\mU$ is generated by $A \cup \{\me_n\:
n\in\N\}$ it follows that $\mU_0$ is dense in $\mU$.
  \proofend

  Recalling the maps $\b_n$ constructed in \lcite{\DefineBeta} we
have:

  \state Proposition
  \label \CommutingDiagram
  For each $n\in\N$ the diagram
  $$
  \matrix{ \Lin(M_n) & \map{\Phi_n} & \Lin(M_{n+1})\cr\cr
           {\scriptstyle \b_n} \downarrow \quad &&
           \qquad \downarrow {\scriptstyle \b_{n+1}} \cr\cr
           \Lin(M_{n+1}) & \map{\Phi_{n+1}} & \Lin(M_{n+2})}
  $$
  commutes and hence there exists a unique injective *-endomorphism
  $$
  \beta : B \to B
  $$
  of the inductive limit C*-algebra $B$, which coincides with $\b_n$ on
each $\Lin(M_n)$.
Moreover
  \izitem
  \zitem $\b(\me_n) = \me_{n+1}$, for all $n\in \N$,
  \zitem $\b(a) = \a(a)\me_1 = \me_1\a(a)$, for all $a\in A$,
  \zitem $\mU$ is invariant under $\b$, and
  \zitem $\b(\mU) = \me_1\mU \me_1$.

  \proof
  Given $T\in \Lin(M_n)$ we have for all $a\in A$ that
  $$
  \Phi_{n+1}(\b_n(T))a =
  \b_n(T)a =
  \ma_n T \mtr_n(a) =
  \ma_n T \Tr(a).
  $$
  On the other hand
  $$
  \b_{n+1}(\Phi_n(T))a =
  \ma_{n+1} \Phi_n(T) \mtr_{n+1}(a) =
  \ma_{n+1} \Phi_n(T) \Tr(a) =
  \ma_{n+1} T \Tr(a).
  $$
  By checking first on the dense set $A\subseteq M_n$
  it is easy to see that the following diagram commutes
  $$
  \matrix{ M_n & \map{\ma_n} & M_{n+1}\cr
           \downarrow && \downarrow  \cr
           M_{n+1} & \map{\ma_{n+1}} & M_{n+2}}
  $$
  where the vertical arrows are the standard embeddings.
In other words $\ma_n(\xi)=\ma_{n+1}(\xi)$ for all $\xi\in M_n$.  If
we now plug $\xi:=T \Tr(a)$ in this identity we conclude that
  $\Phi_{n+1}(\b_n(T)) = \b_{n+1}(\Phi_n(T))$
  as desired.
  Considering (i) we have
  $$
  \b(\me_n) = \b_n(\me_n) = \ma_n \me_n \mtr_n = \me_{n+1},
  $$
  by \lcite{\VariousFormulas.iii}.  Given $a\in A = \Lin(M_0)$, and
$b\in A \subseteq M_1$, we have
  $$
  \b(a)(b) = \b_0(a)(b) = \ma_0 a \mtr_0 (b) =
  \a( a \Tr(b)) =
  \a(a) \a(\Tr(b)) =
  \a(a) E(b) =
  \a(a) \me_1 (b),
  $$
  so we see that $\b(a)= \a(a) \me_1$.  Taking adjoints 
we also have that $\b(a)= \me_1\a(a)$, hence proving (ii).
  It is clear that (iii) follows from (i) and (ii).  In order to prove
(iv) observe that by (ii) we have $\b(1)=\me_1$ so it must be that
$\b(\mU)\subseteq \me_1 \mU \me_1$.

  To prove the reverse inclusion it suffices, by
\lcite{\AlgebraicStructureofF.iii}, to show that for all $a,b\in A$
and $n\in \N$ one has that $\me_1(a \me_n b) \me_1\in\b(\mU)$.  Assuming
initially that $n\geq 1$ notice that
  $$
  \me_1(a \me_n b)\me_1  =
  (\me_1 a \me_1) \me_n (\me_1 b \me_1) =
  E(a) \me_1 \me_n \me_1 E(b) =
  \a(\Tr(a)) \me_1 \me_n \me_1 \a(\Tr(b)) =
  \b\big(\Tr(a) \me_{n-1} \Tr(b)\big).
  $$
  On the other hand if $n=0$ we have
  $$
  \me_1(ab)\me_1  = E(ab)\me_1 = \a(\Tr(ab))\me_1 = \b(\Tr(ab)).
  \proofend
  $$

  Recall (\cite{\CuntzTwo}, \cite{\Stacey}, \cite{\Murphy},
\scite{\Endo}{4.4 and 4.7}) that, since $\b(\mU)$ is a hereditary
subalgebra of $\mU$ by \lcite{\CommutingDiagram.iv}, the crossed
product $\MCP$ is the universal unital C*-algebra generated by a
copy of $\mU$ and an isometry $\mS$ subject to the relation that
  $$
  \mS x\mS^* = \b(x) \for x\in \mU.
  $$
  It is well known  that $\mU$ embeds injectively in $\MCP$
(see the remark after Definition 4.4 in
\cite{\Endo}, \scite{\Stacey}{Section 2}, or \scite{\Murphy}{Section
2}), so we will view $\mU$ as a subalgebra of $\MCP$.  Since $A$ is a
subalgebra of $\mU$ we also have that $A\subseteq \MCP$.

  \state Proposition
  \label \ClosedProblem
  There exists a *-homomorphism $\phi: \CP \to \MCP$ such that
$\phi(\cS ) = \mS$, and $\phi(q(a))=a$, for all $a\in A$, where $q$ is
the canonical  quotient map from $\TCP$ to $\CP$.

  \proof
  It will be useful to keep in mind that
  $\me_1 = \b(1) = \mS\mS^*$, and hence that $\me_1 \mS = \mS.$
  Considering the natural inclusion of $A$ in $\mU$ notice that for
all $a\in A$ one has
  $$
  \mS a =
  \mS a \mS^* \mS =
  \b(a) \mS =
  \a(a) \me_1 \mS =
  \a(a) \mS.
  $$
  Also
  $$
  \mS^*a \mS =
  \mS^*\me_1a \me_1\mS =
  \mS^*E(a) \me_1\mS =
  \mS^*E(a)\mS =
  \mS^*\a(\Tr(a))\mS =
  \mS^*\mS \Tr(a)=
  \Tr(a).
  $$
  It follows from the universal property of $\TCP$
that there exists a *-homomorphism $\phi:\TCP \to \MCP$
which is the identity on $A$ and such that $\phi(\tS )=\mS$.

Now let $(a,k)$ be a redundancy in $\TCP$.  We claim that $a=\phi(k)$.
In order to see this note that since $k\in \overline{A\tS \tS^*A}$ one
has that $\phi(k) \in \overline{A\mS \mS^*A} =
\overline{A\me_1A}$.

  For all $b\in A$ it is assumed that $ab\tS =kb\tS $ so
that $ab\tS \tS^*=kb\tS \tS^*$ and hence $ab\me_1=\phi(k)b\me_1$.
  Observe that all terms occurring in this last  identity lie in the
algebra generated by $A$ and $\me_1$, which consists of
operators on the Hilbert module $M_1$.  In particular, considering 1
as an element of $M_1$ we have that
  $$
  ab = ab\me_1(1) = \phi(k)b\me_1(1) = \phi(k)b.
  $$
  If follows that $a$ and $\phi(k)$ coincide on $A$, which is a dense
subspace of $M_1$, and hence that $a=\phi(k)$ as claimed.
  Therefore $\phi$ vanishes on the ideal generated by the differences
$a-k$ for all redundancies $(a,k)$ and hence factors through the
quotient yielding the desired map.
  \proofend

  We may now answer a question raised in \cite{\Endo}:

  \state Theorem
  Let $\a$ be an injective endomorphism of a unital C*-algebra $A$
with $\a(1)=1$ and let $\Tr$ be a transfer operator of the form
$\Tr=\a\inv\compos E$, where $E$ is a conditional expectation from $A$
to the range of $\a$.  Then the natural map from $A$ to $\CP$ is
injective.

  \proof
  Given that $\phi(q(a))=a$, we see that $q$ is injective when
restricted to $A$.
  \proofend

 From now on we are therefore allowed to identify $A$ and its image
$q(A)$ within $\CP$.

  \section{Consequences of Watatani's work}
  Our next main goal will be to prove that the map $\phi$ of
\lcite{\ClosedProblem} is in fact an isomorphism.  But before we are
able to attack this question we need to do some more work.
We begin by introducing some notation.  Recalling that
$\tS $ is the standard isometry in $\TCP$
we will let
  $$
  \te_n = \tS^n\tS^{*n},
  $$
  and $\tK_n$ be the closed linear span of $A \te_n $A, so
that each $\tK_n$ is a sub-C*-algebra of $\TCP$, as well as an
$A$--bimodule.
Some elementary properties of the $\te_n$ and the $\tK_n$ are recorded
in the next:

  \state Proposition
  \label \ElementaryPropertiesOfKn
  For all $n,m\in \N$ with $n\leq m$ one has that
  \izitem
  \zitem $\te_m \leq \te_n$,
  \zitem $\overline{\tK_n \tK_m} = \overline{\tK_m \tK_n} = \tK_m$,
  \zitem $\te_n a \te_n = \E_n(a) \te_n =  \te_n \E_n(a)$, for all $a\in A$,
  \zitem $\tS^*\tK_{n+1}\tS  \subseteq \tK_n$, and
  \zitem $\tS \tK_n \subseteq \tK_{n+1}\tS $.

  \proof
  The first point is trivial and hence we omit it.  Regarding (iii)  we
have by induction that
  $$
  \te_{n+1}a  \te_{n+1} =
  \tS  \te_n\tS^* a \tS  \te_n \tS^* =
  \tS  \te_n \Tr(a) \te_n \tS^* =
  \tS  \E_n(\Tr(a))\te_n \tS^* \&=
  \a \E_n \Tr (a) \te_{n+1} =
  \a \E_n \a\inv E (a) \te_{n+1} =
  \E_{n+1} (a) \te_{n+1}.
  $$
  This proves that $\te_na \te_n = \E_n (a) \te_n$ and hence also that
$\te_na \te_n = \te_n \E_n (a)$ by taking adjoints.
  Speaking of (ii) notice that
  $$
  (A\te_n A)  (A\te_m A) =
  A(\te_n A\te_n) \te_m A \subseteq
  A\E_n(A)\te_n \te_m A =
  A\E_n(A) \te_m A \subseteq
  \tK_m,
  $$
  proving that $\tK_n \tK_m \subseteq \tK_m$ so that $\overline{\tK_n
\tK_m} \subseteq \tK_m$.  In order to prove the reverse inclusion notice
that for all $a,b\in A$ one has
  $$
  a \te_m b =   (a \te_n 1) (1 \te_m b) \in \tK_n \tK_m.
  $$
  The equality
  $\overline{\tK_m \tK_n} = \tK_m$ follows by taking adjoints.
  Turning now to the proof of   (iv) we have
for all
$a,b\in A$ that
  $$
  \tS^* a \te_{n+1} b \tS  =
  \tS^* a \tS  \te_n \tS^* b \tS  =
  \Tr(a) \te_n \Tr(b)
  \in \tK_n.
  $$
  As for (v) notice that
  $$
  \tS  a\te_n b =
  \a(a)\tS \te_n \tS^*\tS  b =
  \a(a)\te_{n+1}\a(b)\tS  \in
  \tK_{n+1} \tS .
  \proofend
  $$

It is now  convenient to have in mind the sequence
of *-homomorphisms
  $$
  \TCP \map q \CP \map \phi \MCP.
  $$

  \state Lemma
  \label \PreIsomorphismOfKs
  For each $n\in\N$  let
  $\cK_n = q(\tK_n)$,  $\mK_n = \phi(\cK_n)$, and $\ce_n = q(\te_n)$. Then 
  \izitem
  \zitem $\phi(\ce_n)=\me_n$, 
  \zitem $\cK_n= \overline{A \ce_n A}$, and
  \zitem $\mK_n= \overline{A \me_n A}$.

  \proof To prove (i) notice that
  $$
  \phi(\ce_n) =
  \phi(q(\te_n)) =
  \phi(q(\tS^n\tS^{*n})) =
  \mS^n\mS^{*n} =
  \b^n(1) =
  \me_n.
  $$
  As for (ii)
  $
  \cK_n=
  q(\overline{A \te_n A}) = 
  \overline{A \ce_n A}.
  $
  Finally
  $
  \mK_n =
  \phi(\cK_n) =
  \phi(\overline{A \ce_n A}) =
  \overline{A \me_n A}.
  $
  \proofend

  Observe that by \lcite{\PreIsomorphismOfKs.iii} $\mK_n$ is precisely
the \stress{reduced C*-basic construction} \scite{\Watatani}{2.1.2}
relative to the conditional expectation $\E_n: A \to \R_n$.

  In trying to prove that the map $\phi$ of \lcite{\ClosedProblem} is
injective a crucial step will be taken by the following important
consequence of \cite{\Watatani}.

  \state Proposition
  \label \IsomorphismOfKs
  The following *-homomorphisms are in fact
*-isomorphisms
  \izitem
  \zitem $q : \tK_n \to \cK_n$,
  \zitem $\phi : \cK_n \to \mK_n$,
  \zitem $\phi\compos q : \tK_n \to \mK_n$.

  \proof 
  By \scite{\Watatani}{2.2.9} we have that $\mK_n$ is canonically
isomorphic to the \stress{unreduced C*-basic construction}
relative to $\E_n$ and thus 
possesses the universal property described in
\scite{\Watatani}{2.2.7}.

Supposing that $\TCP$ is faithfully represented on a Hilbert space
${\cal H}$ observe that by \lcite{\ElementaryPropertiesOfKn.iii} the
triple
  $(id_A,\te_n,{\cal H})$
  is a \stress{covariant representation} of the conditional
expectation $\E_n$, according to Definition 2.2.6 in \cite{\Watatani}.
It follows that there exists a *-representation $\rho$ of $\mK_n$ on
${\cal H}$ such that $\rho(a\me_nb)=a\te_nb$, for all $a,b\in A$.  Since
$\phi\compos q$ maps $a\te_nb$ to $a\me_nb$ we see that $\phi\compos q$
and $\rho$ are each others inverse, hence proving (iii).  This implies
that $q$ is injective on $\tK_n$ and since $q$ is obviously also
surjective (i) is proven.  Clearly (ii) follows from (i) and (iii).
\proofend

The following elementary properties should also be noted:

  \state Proposition
  \label \MultiModule
  Let $n,m\in\N$ with $n\leq m$.  Then
  \izitem
  \zitem $\overline{\cK_n \cK_m} = \overline{\cK_m \cK_n} = \cK_m$,
  \zitem $\overline{\mK_n \mK_m} = \overline{\mK_m \mK_n} = \mK_m$,
  \zitem Denote by $\psi_n : X_n \to Y_n$ any one of the
isomorphisms in \lcite{\IsomorphismOfKs.i--iii}.  Then
for every $x_n\in X_n$ and $x_m\in X_m$ one has that
  $\psi_m(x_n x_m) = \psi_n(x_n)\psi_m(x_m)$ and
  $\psi_m(x_m x_n) = \psi_m(x_m)\psi_n(x_n)$.

Taking $n=0$ in \lcite{\MultiModule.iii}, in which case
$\tK_n=\cK_n=\mK_n=A$, we see that all of the isomorphisms in
\lcite{\IsomorphismOfKs.i--iii} are also $A$--bimodule maps.

  \section{Higher Order Redundancies}
  We now wish to study a generalization of the notion of redundancy.
For this we need the following fact for which we have found no
reference in the literature.

  \state Lemma
  \label \UsefulTechnical
  Let $B$ and $J$ be closed *-subalgebras of some C*-algebra $C$ such
that $\overline{JB}=J$.  If $x\in J$ and $xB\subseteq B$ then $x\in
B$.

  \proof
  Viewing $J$ as a right Banach $B$--module we have by the
Cohen-Hewitt factorization theorem \scite{HR}{32.22} that $x=ya$ for
some $y\in J$ and $a\in B$.  Choosing an approximate unit $\{u_i\}_i$
for $B$ we have that
  $$
  x = ya = \lim_{i\to\infty} yau_i= \lim_{i\to\infty} xu_i \in
\overline{xB} \subseteq B.
  \proofend
  $$

  \definition
  \label \DefineNRedundancy
  Let $n\geq1$ be an integer.  A \stress{redundancy of order $n$}, or
an \stress{$n$--redundancy}, is a finite sequence $(a_0,a_1,\ldots,a_n)\in
\prod_{i=0}^n\tK_i$ such that $\sum_{i=0}^n a_ix=0$, for all $x\in
\tK_n$.

Up to a minus sign the above notion generalizes the notion of
redundancy introduced in \cite{\Endo}.  In fact it is easy to see that
the pair $(a,k)$ is a redundancy according to \cite{\Endo} if and only if
$(a,-k)$ is a $1$--redundancy.

We now come to a main technical result:

  \state Proposition
  \label \BigTech
  Let $n\geq1$  and let $(a_0,a_1,\ldots,a_n)$ be a redundancy of order $n$.  Then
$\sum_{i=0}^n q(a_i) = 0$.

  \proof
  We proceed by induction observing that the case $n=1$ follows
easily from
the observation already made that $1$--redundancies are simply
redundancies.
  So let $n>1$ and let $(a_0,a_1,\ldots,a_n)$ be an $n$--redundancy.
Given $b$ in $A$ let $a'_i= \tS^*b^*a_ib\tS $ for all $i=0,1,\ldots,n$ and
observe that by \lcite{\ElementaryPropertiesOfKn.v} one has that
  $$
  \(\sum_{i=0}^n a'_i\) \tK_{n-1} =
  \tS^*b^*\(\sum_{i=0}^n a_i\)b \tS  \tK_{n-1} \subseteq
  \tS^*b^*\(\sum_{i=0}^n a_i\)b \tK_n \tS   = \{0\}.
  $$
  Since  $a_0' = \Tr(b^*a_0b)\in A$, and
  $a_i'\in \tK_{i-1}$ for $i\geq 1$, by \lcite{\ElementaryPropertiesOfKn.iv}, we have
that $(a'_0+a'_1,a'_2,\ldots,a'_n)$ is a redundancy of order $n-1$.  So
  $\sum_{i=0}^n q(a'_i)=0$
  by induction.
  Equivalently $\sum_{i=0}^n a'_i \in\Ker(q)$.  Assume first that
$\sum_{i=0}^n a_i$ is positive.  We then have that
  $$
  \Ker(q)\ \ni\ \sum_{i=0}^n a'_i =
  \tS^*b^*\(\sum_{i=0}^n a_i\)b\tS  =
  \tS^*b^*\(\sum_{i=0}^n a_i\)^{1/2}\(\sum_{i=0}^n a_i\)^{1/2}b\tS .
  $$
  It follows that $\(\sum_{i=0}^n a_i\)^{1/2}b\tS \in\Ker(q)$ and hence also
$\(\sum_{i=0}^n a_i\)b\tS \in\Ker(q)$.
Multiplying this on the right by $\tS^{n-2}(\tS^*)^{n-1}c$, for $c\in A$,
we see that
  $$
  \(\sum_{i=0}^n a_i\) b\te_{n-1}c \ \in\ \Ker(q) \for b,c\in A,
  $$
  and hence that $\(\sum_{i=0}^n a_i\)\tK_{n-1}\subseteq\Ker(q)$.  For
all $y\in q(\tK_{n-1})= \cK_{n-1}$ it follows that
  $$
  q(a_n) y = - \sum_{i=0}^{n-1} q(a_i) y,
  \eqno{(\seqnumbering)}
  \label \AnEquation
  $$
  from where we deduce that  $q(a_n) \cK_{n-1} \subseteq \cK_{n-1}$.   By
\lcite{\UsefulTechnical} with $B=\cK_{n-1}$ and $J=\cK_n$ we have that
  $q(a_n)\in\cK_{n-1}$ and hence there exists $b_n\in\tK_{n-1}$ such that
$q(b_n)=q(a_n)$.
  Observe that for every $x\in \tK_{n-1}$ we then have
  $$
  q\(\sum_{i=0}^{n-1} a_i x + b_nx\) =
  \sum_{i=0}^{n-1} q(a_i) q(x)  + q(a_n)q(x) =
  0,
  $$
  by \lcite{\AnEquation}.
  Observing that the term within the big parenthesis above lies in
$\tK_{n-1}$, we have   by \lcite{\IsomorphismOfKs.i} that
  $$
  \sum_{i=0}^{n-1} a_i x + b_nx = 0 \for x\in\tK_{n-1},
  $$
  and hence
  $(a_0,a_1,\ldots,a_{n-2}, a_{n-1} + b_n)$
  is a redundancy of order $n-1$.
Once again by the induction
hypothesis it follows that
  $$
  0 =
  \sum_{i=0}^{n-1} q(a_i) + q(b_n)  =
  \sum_{i=0}^n q(a_i).
  $$

Without assuming that $\sum_{i=0}^n a_i$ be
positive one can expand the expression
  $\(\sum_{i=0}^n a_i\)^*\(\sum_{i=0}^n a_i\)$
  and rearrange its terms   in order to form a redundancy
  $(b_0,b_1,\ldots,b_n)$ such that
  $$\sum_{i=0}^n b_i = \(\sum_{i=0}^n a_i\)^*\(\sum_{i=0}^n a_i\)$$
  and the conclusion will follow easily.
  \proofend

  \state Theorem
  \label \PhiIsIsomorphism
  The map $\phi: \CP \to \MCP$ of \lcite{\ClosedProblem} is an isomorphism.

  \proof
  We begin by proving that $\phi$ is surjective.
  Since $\phi$ is the identity on $A$ and since $\phi(\cS^n\cS^{*n}) =
\mS^n1\mS^{*n} = \b^n(1) = \me_n$ we have that $\mU$, which is generated
by $A \cup \{\me_n\: n\in\N\}$, is contained in the range of $\phi$.  On
the other hand $\phi(\cS)=\mS$ and $\MCP$ is generated by
$\mU\cup\{\mS\}$.  So we see that $\phi$ is indeed surjective.

Using the universal property of $\MCP$ it is easy to see that there
exists a circle action $\mgamma$  on $\MCP$ such that
  $$
  \mgamma_z(\mS) = z\mS\and
  \mgamma_z(f) = f
  \for f\in \mU
  \for z\in S^1.
  $$
  Since $\phi$ is clearly covariant with respect to $\cgamma$ and
$\mgamma$, if we prove that $\phi$ is injective on the fixed point
subalgebra for $\cgamma$, which we denote by $F$, then by
\scite{\newpim}{2.9} we would have proven that $\phi$ is injective.
Recall from \lcite{\FixedPointAlgebras} that 
  $
  F =
  \overline{\sl span}\{a\cS^n\cS^{*n}b: a,b\in A,\ n\in \N\}.
  $
  If we further observe that
  $
  \cS^n\cS^{*n} =
  q(\tS^n\tS^{*n}) =
  q(\te_n) =
  \ce_n
  $
  we see that
  $$
  F = \overline{\sum_{n\in\N} \cK_n}.
  $$
  In order to prove that $\phi$ is injective it is thus enough to show it
to be injective, and hence isometric,  on each subalgebra of the form
  $$
  F_n = \overline{\sum_{0\leq i\leq n} \cK_i},
  $$
  for $n\in \N$.  Applying \scite{\Ped}{1.5.8} repeatedly it is easy
to see that $\sum_{0\leq i\leq n} \cK_i$ is closed so that, in fact,
  $F_n = \sum_{0\leq i\leq n} \cK_i$.

  Let $a\in F_n$ be such that $\phi(a)=0$ and write $a = \sum_{i=0}^n
a_i$, with $a_i\in\cK_i$.   For every $i=0,\ldots,n$ choose
$b_i\in\tK_i$ with $q(b_i)=a_i$ (such  $b_i$ exists and is unique
by \lcite{\IsomorphismOfKs.i}).  We now claim that $(b_0,b_1,\ldots,b_n)$ is a
redundancy of order $n$.  In fact, for every $x\in\tK_n$ one has that
  $$
  \phi\compos q\(\sum_{i=0}^n b_ix\) =
  \phi\(\sum_{i=0}^n a_i\)\phi(q(x)) =
  \phi(a)\phi(q(x)) =0.
  $$
  By \lcite{\IsomorphismOfKs.iii} we have that
  $\sum_{i=0}^n b_ix=0$ hence proving our claim.  Employing
\lcite{\BigTech} we then conclude that
  $
  a = \sum_{i=0}^n q(b_i) = 0.
  $
  \proofend

  \section{A bit of cohomology}
  We will set this section aside to list certain elementary definitions
and facts about an ingredient of cohomological flavor which will be
recurrent in our development.

  \definition
  \label \CohomologyDefinition
  Given $a\in A$ and $n\in \N$ we will let
  $$
  a^{[n]}:= a\; \a(a) \cdots \a^{n-1}(a)
  $$
  with the convention that $a^{[0]} = 1$.

  It is elementary to prove that:
  
  \state Proposition
  \label \CohomologyBoundary
  For all $a\in A$ and $n,m\in\N$ one has that
  $a^{[n+m]} = a^{[n]} \a^n( a^{[m]})$,
  
  Although $\Zenter(A)$  is not necessarily invariant by $\a$ observe
that given $a,b\in\Zenter(A)$ and $n,m\in\N$ one has that $\a^n(a)$
and $\a^m(b)$ commute.  In fact, supposing without loss of generality
that $n\leq m$, observe that
  $$
  \a^n(a) \a^m(b) =
  \a^n\big(a \a^{m-n}(b)\big) =
  \a^n\big(\a^{m-n}(b)a\big) =
  \a^m(b)\a^n(a).
  $$

  This shows that, when $a\in\Zenter(A)$, the order of the factors in the
definition of $a^{[n]}$ above is irrelevant.   It is also easy to
conclude that:

  \state Proposition
  \label \CohomologyMiscelanea
  For all $a,b\in\Zenter(A)$ and $n\in\N$ one has that 
  \izitem
  \zitem $a^{[n]}b^{[n]}=(ab)^{[n]}$,
  \zitem if $a$ is invertible then $(a^{[n]})\inv =(a\inv)^{[n]}$,
  \zitem if $a$ is a self-adjoint and 
  $0\leq a\leq c$, where $c\in\Real$, then $0\leq a^{[n]}\leq c^n$.
  %  \zitem if $a\geq0$ then $a^{[n]}\geq 0$.

  \definition
  From now on we will denote by $\comm{\R_n}{A}$ the commutant of
$\R_n$ in $A$, that is, the set of elements in $A$ which commute with
$\R_n$.

  Since the $\R_n$ are decreasing it is clear that the
$\comm{\R_n}{A}$ are increasing.  Moreover, if $a\in\Zenter(A)$ it is
clear that $\a^k(a) \in \Zenter(\R_k) \subseteq \comm{\R_k}{A}$.  From
this one immediately has:

  \state Proposition
  \label \InTheRightComm
  For all $a\in\Zenter(A)$ and $n\geq 1$ one has that 
  $a^{[n]}\in \comm{\R_{n-1}}{A}$.

  In connection with our transfer operator $\Tr$ we will later need
the following fact:
  
  \state Proposition
  \label \FormuLoca
  For every invertible element $\lambda\in\Zenter(A)$, $a\in A$,  and $n,m,p\in\N$ one has
  \izitem
  \zitem
  $
  \lambda^{-[m+p]} \a^m\Tr^n( \lambda^{[n+p]}a )= 
  \lambda^{-[m]} \a^m\Tr^n( \lambda^{[n]}a ),  
  $
  \zitem
  $
  \a^m\Tr^n(a \lambda^{[n+p]} )\lambda^{-[m+p]}= 
  \a^m\Tr^n(a \lambda^{[n]} )   \lambda^{-[m]},  
  $
  \medskip\noindent
  where by $\lambda^{-[m]}$ we mean $(\lambda\inv)^{[m]}$.

  \proof Observing that for all $x,y\in A$ we have 
  $$
  \a^m\Tr^n(\a^n(x)y) =
  \a^m(x\Tr^n(y)) = 
  \a^m(x)\a^m\Tr^n(y),
  $$
  and 
using \lcite{\CohomologyBoundary} to compute
  $\lambda^{-[m+p]}$ and   $\lambda^{-[n+p]}$  we have
  $$
  \lambda^{-[m+p]} 
  \a^m\Tr^n( \lambda^{[n+p]}a )= 
  \lambda^{-[m]} \a^m(\lambda^{-[p]}) 
  \a^m\Tr^n(\lambda^{[n]} \a^n(\lambda^{[p]})a ) =
  \lambda^{-[m]} 
  \a^m\Tr^n(\lambda^{[n]} a ).
  $$
  This proves (i) and the proof of (ii) is similar.
  \proofend

  \section{Finite Index automorphisms}
  \label \FiniteIndexSection
  Given a C*-algebra $B$ and a closed *-subalgebra $C$ recall from
\scite{\Watatani}{1.2.2 and 2.1.6} that a conditional expectation $E:
B \to C$ is said to be of \stress{index-finite type} if there exists a
\stress{quasi-basis} for $E$, i.e.~a finite sequence
$\{u_1,\ldots,u_m\}\subseteq B$ such that
  $$
  a=\sum_{i=1}^m u_i E(u_i^* a)\for a\in B.
  $$
  In this case one defines  the \stress{index} of $E$ by
  $$
  \ind(E) = \sum_{i=1}^m u_iu_i^*.
  $$
  It is well known that $\ind(E)$
  does not depend on the choice of the $u_i$'s, 
  that it belongs to the center of $B$ \scite{\Watatani}{1.2.8} 
  and
  is invertible \scite{\Watatani}{2.3.1}.

  \definition We shall say that a pair $(\a,E)$ is a \stress{finite index
endomorphism} of the C*-algebra $A$ if $\a$ is a *-endomorphism of $A$
and $E$ is a conditional expectation of index-finite type from $A$ to the
range of $\a$.

Throughout this section we will fix a finite index endomorphism
$(\a,E)$ and a quasi-basis $\{u_1,\ldots,u_m\}$ for $E$.  As before we
will let $\Tr$ be the transfer operator given by $\Tr = \a\inv\compos
E$.

  \state Proposition
  \label \EnFunctionOfEnPlusOne
  For every $n\in\N$ one has that
  $\displaystyle \me_n=\sum_{i=1}^m \a^n(u_i) \me_{n+1} \a^n(u_i^*)$.

  \proof
  Observe
that for all $a\in A \subseteq M_1$ one has that
  $$
  \sum_{i=1}^m u_i \me_1 u_i^* (a) =
  \sum_{i=1}^m u_i E( u_i^*a) =
  a,
  $$
  so that
  $1=\sum_{i=1}^m u_i \me_1 u_i^*$, hence proving the statement for
$n=0$.
  Assuming that $n\geq1$ apply the endomorphism $\b$ of
\lcite{\CommutingDiagram} to both sides of the expression
  $\me_{n-1}=\sum_{i=1}^m \a^{n-1}(u_i) \me_n \a^{n-1}(u_i^*)$
  to conclude that
  $$
  \me_n = \b(\me_{n-1}) =
  \sum_{i=1}^m \b(\a^{n-1}(u_i)) \me_{n+1} \b(\a^{n-1}(u_i^*)) \&=
  \sum_{i=1}^m \a^n(u_i)\me_1 \me_{n+1} \me_1\a^n(u_i^*) =
  \sum_{i=1}^m \a^n(u_i) \me_{n+1} \a^n(u_i^*).
  \proofend
  $$

  As a consequence we have:

  \state Proposition
  \label \KnIncreasing
  For every $n\in\N$ one has that $\mK_n \subseteq \mK_{n+1}$.

  \proof
  Given $a,b\in A$ observe that by
\lcite{\EnFunctionOfEnPlusOne} one has that
  $$
  a\me_nb =
  \sum_{i=1}^m a\; \a^n(u_i) \me_{n+1} \a^n(u_i^*)\;b,
  $$
  and hence $A \me_n A \subseteq A \me_{n+1} A$.  The conclusion then
follows from \lcite{\PreIsomorphismOfKs.iii}.
  \proofend

  It is our next major goal to define a conditional expectation $G:\CP
\to A$.  In order to do this we will first define conditional expectations
  $G_n: \mK_n \to A$.  Considering that any such conditional
expectation is an $A$--bimodule map we will begin with a brief study
of $A$--bimodule maps.  So let $f:\mK_n\to A$ be any $A$--bimodule map.
Observing that $\me_n$ commutes with $\R_n$ by \lcite{\AlgebraicStructureofF.ii}
notice that for all $x\in
\R_n$ we have that
  $$
  xf(\me_n) =
  f(x\me_n) =
  f(\me_nx) =
  f(\me_n)x,
  $$
  so $f(\me_n)\in \comm{\R_n}{A}$.

  \state Proposition
  \label \TensorMap
  For every $\lambda\in \comm{\R_n}{A}$ there is one and only one
$A$--bimodule map $f:\mK_n\to A$ such that
  $$
  f(a\me_nb)=a\lambda b \for a,b\in A.
  $$
  If
$\lambda\geq0$ then $f$ is a positive map and vice-versa.

  \proof
  Recall from \scite{\Watatani}{2.2.2} that $\mK_n$, being the C*-basic
construction relative to the conditional expectation $\E_n: A \to
\R_n$, is isomorphic to the algebraic tensor product (no completion)
$\tensor$ under the map
  $$
  \psi:  a\* b \in \tensor \longmapsto a \me_n b \in \mK_n.
  $$
  Given $\lambda\in \comm{\R_n}{A}$ we have that
the map
  $$
  (a,b)\in A\x A \mapsto a\lambda b\in A
  $$
  is $\R_n$--balanced, which in turn defines a linear map $\tilde f:\tensor\to
A$ such that $\tilde f(a\*b) = a\lambda b$, for all $a$ and $b$ in
$A$.  Composing this with the inverse of the isomorphism $\psi$
defined above yields a linear map $f:\mK_n\to A$ such that
$f(a\me_nb) = a\lambda b$.  It is now easy to see that $f$ is an $A$--bimodule
map and that it is uniquely determined in terms of $\lambda$.

  If $f$ is positive it is clear that $\lambda = f(\me_n)\geq0$.
Conversely suppose that $\lambda\geq0$.  Given $x\in\mK_n$ 
of the form 
  $
  x = \sum_{i=1}^m a_i \me_n b_i
  $
  notice that
  $$
  x^*x  =
  \sum_{i,j=1}^m b_i^* \me_n a_i^* a_j \me_n b_j =
  \sum_{i,j=1}^m b_i^* \E_n(a_i^* a_j)\me_n b_j,
  $$
  by \lcite{\AlgebraicStructureofF.ii}.
  Since conditional expectations are completely positive by
\scite{\Takesaki}{IV.3.4}, we have that $\{\E_n(a_i^* a_j)\}_{ij}$ is
a positive matrix and hence there exists an $m\times m$ matrix $c =
\{c_{ij}\}_{ij}$ over $\R_n$ such that
  $$
  \E_n(a_i^* a_j) =
  \sum_{k=1}^m c^*_{ki}c_{kj}
  \for i,j=1,\ldots,m.
  $$
  Therefore 
  $$
  x^*x  =
  \sum_{i,j,k=1}^m b_i^* c^*_{ki} \me_n  c_{kj}b_j =
  \sum_{k=1}^m d_k^* \me_n d_k,
  $$
  where 
  $d_k = \sum_{j=1}^m c_{kj}b_j$.  If follows that
  $$
  f(x^*x)  =
  \sum_{k=1}^m d_k^* \lambda d_k \geq0.
  \proofend
  $$

  As mentioned above $\ind(E)\in \Zenter(A)$ so we have  by
\lcite{\InTheRightComm} that
  $$
  \iota_n := \big(\ind(E)\big)^{[n]} \in
  \comm{\R_{n-1}}{A} \subseteq
  \comm{\R_n}{A}
  \for n\geq1.
  \eqno{(\seqnumbering)}
  \label \DefineIn
  $$

  \state Proposition
  \label \GnRestrito
  For each $n\in\N$   let 
  $
  G_n : \mK_n \to A
  $
  be the unique $A$--bimodule map 
  such that
  $$
  G_n(a\me_nb)=a\iota_n\inv b\for a,b\in A,
  $$
  given by 
\lcite{\TensorMap}.  Then
  for every $n\in\N$ one has that $G_n$ is a positive
contractive conditional expectation from $\mK_n$ to $A$.  Moreover
the restriction of $G_{n+1}$ to $\mK_n$ coincides with $G_n$.

  \proof
  By \lcite{\EnFunctionOfEnPlusOne} we have that
  $$
  G_{n+1}(\me_n) =
  G_{n+1}\(\sum_{i=1}^m \a^n(u_i) \me_{n+1} \a^n(u_i^*)\) =
  \sum_{i=1}^m \a^n(u_i) \iota_{n+1}\inv \a^n(u_i^*) = \cdots
  $$
  Observing that $\iota_{n+1}\in \comm{\R_n}{A}$ by \lcite{\DefineIn}
we see that the above
equals
  $$
  \cdots =
  \sum_{i=1}^m \a^n(u_i) \a^n(u_i^*) \iota_{n+1}\inv =
  \a^n(\ind(E)) \iota_{n+1}\inv =
  \iota_n\inv =
  G_n(\me_n).
  $$
  This proves that $G_{n+1}(\me_n)= G_n(\me_n)$ from where one easily
deduces that $G_{n+1}|_{\mK_n}= G_n$.  It follows that each $G_n$
coincides with $G_0$ on $\mK_0$.  In other words $G_n$ is the identity
on $A$ and, being an $A$--bimodule map, we see that it is in fact a
conditional expectation onto $A$.  Since $I_n$ is positive by
\lcite{\CohomologyMiscelanea.iii} we have by \lcite{\TensorMap} that
$G_n$ is a positive map.  To conclude observe that a positive
conditional expectation is always contractive.
  \proofend

  \sysstate{Remark}{\rm}
  {We should remark that $G_1$ is precisely the dual conditional
expectation defined in \scite{\Watatani}{2.3.2}.  For $n\geq2$, even
though each $\E_n$ is a conditional expectation of index-finite type
and $\mK_n$ is the C*-basic construction relative to $\E_n$, $G_n$ may
not be the dual conditional expectation if $A$ is non-commutative.
This is due to the fact that $\ind(\E_n)$ may differ from $\iota_n$.
See also \scite{\Watatani}{1.7.1}.}

Proposition \lcite{\GnRestrito} says that the $G_n$ are compatible
with each other and hence may be put together in the following way:

  \state Proposition
  \label \ConditionalExpectOnF
  There exists a conditional expectation $\mF:\mU \to A$ such that
  $$
  \mF(a\me_nb) = a\iota_n\inv b \for a,b\in A \for n\in\N.
  $$
  If $A$ is commutative there is no other conditional expectation from
$\mU$ to $A$.

  \proof From \lcite{\AlgebraicStructureofF.iii} it follows that
  $\mU =
  \overline {\pilar 9 \sum_{n\in\N} A \me_n A} =
  \overline {\pilar 9\sum_{n\in\N} \mK_n}$, but since the $\mK_n$ are
increasing by \lcite{\KnIncreasing}, we see that $\mU$ is in fact the
inductive limit of the $\mK_n$.  The existence of $\mF$ then follows
easily from \lcite{\GnRestrito}.

Since $\mK_n$ is the C*-basic construction relative to $\E_n$ it
follows from \scite{\Watatani}{1.6.4} that there exists a unique conditional
expectation from $\mK_n$ to $A$, under the hypothesis that $A$ is
commutative.  Therefore any conditional expectation $F'$ from $\mU$ to
$A$ must coincide with $G_n$ on each $\mK_n$ and hence $F'=\mF$.
  \proofend

  We now come to one of our main results:

  \state Theorem
  \label \UniqueCondExpInvariant
  Let $A$ be a unital C*-algebra and let $(\a,E)$ be a finite index
endomorphism of $A$ such that $\a$ is injective and preserves the
unit.  Then there exists a conditional expectation $G:\CP\to A$ such
that
  $$
  G(a\cS^n\cS^{*m}b) =
  \delta_{nm}a\iota_n\inv b
  \for a,b\in A \for n,m\in\N,
  $$
  where $\delta$ is the Kronecker symbol.
  If $A$ is commutative any conditional expectation from $\CP$ to $A$
which is invariant under the scalar gauge action $\cgamma$ (see
\lcite{\ScalarGaugeAction}) coincides with $G$.

  \proof
  Identifying $\CP$ and $\MCP$ via the isomorphism $\phi$ of
\lcite{\PhiIsIsomorphism} it is enough to prove the corresponding
result for $\MCP$, with $\mS$ replacing $\cS$, and the action
$\mgamma$
described in the proof of \lcite{\PhiIsIsomorphism} replacing
$\cgamma$.
  Consider the operator
  $\check P$ on $\MCP$ given by
  $$
  \check P(a) = \int_{S^1} \mgamma_z (a)\, dz \for a\in\MCP.
  $$
  As mentioned in the proof of \lcite{\PreFixedPointAlgebras} $\check
P$ is a conditional expectation onto the fixed point algebra for
$\mgamma$.   Considering that $\MCP$ is the closed linear span of the
set
  $
  \{a\mS^n\mS^{*m}b: a,b\in A,\ n,m\in \N\}
  $
  one may use \lcite{\PreFixedPointAlgebras} in order to prove that
the fixed point algebra for $\mgamma$ coincides with $\mU$.  The
composition
  $G:= \mF\compos \check P$
  is therefore the conditional expectation sought.  Now suppose that
$A$ is commutative and let $G'$ be any conditional expectation from
$\MCP$ to $A$.  By \lcite{\ConditionalExpectOnF} we have that $G'|_\mU
= \mF$.  If $G'$ is  moreover invariant
under $\mgamma$ we have for all $a\in\MCP$ that
  $$
  G'(a)  =
  G'\(\int_{S^1} \mgamma_z (a)\, dz \) =
  G'(\check P(a)) =
  \mF(\check P(a)) =
  G(a).
  \proofend
  $$

  As an immediate consequence we have:

  \state Corollary
  There exists a conditional expectation $\tG:\TCP\to A$ such that
  $$
  \tG(a\tS^n\tS^{*m}b) =
  \delta_{nm}a\iota_n\inv b
  \for a,b\in A \for n,m\in\N.
  $$

  \proof
  It is enough to put $\tG=G\compos q$.
  \proofend

  \section{KMS states}
  \label \KMSSection
  Throughout this section and until further notice we will assume
the following:

  \sysstate{Standing Hypotheses}{\rm}{\label \StandingHyp 
  \izitem
  \zitem $A$ is a unital C*-algebra,
  \zitem $\a$ is an injective endomorphism of $A$ such that $\a(1)=1$,
  \zitem $E$ is a conditional expectation from $A$ to the range of $\a$,
  \zitem $E$ is of index-finite type,
  \zitem $\Tr$ is the transfer operator given by $\Tr = \a\inv\compos
E$,
  \zitem $h$ is a fixed self-adjoint element in the center of $A$ such
that $h \geq cI$ for some real number $c>0$ (later we will actually require
that $c>1$),
  \zitem $\tsigma$ and $\csigma$ will denote the gauge actions
referred to in \lcite{\DefineGaugeAction} as $\tsigma^h$ and
$\csigma^h$, respectively.}

  \bigskip
  The purpose of this section will be to study the KMS states on
$\TCP$ and $\CP$ relative to $\tsigma$ and $\csigma$.  Whenever we say
that a state is a KMS state on $\TCP$ (resp.~$\CP$) it will be with
respect to $\tsigma$ (resp.~$\csigma$).

Observe that the canonical quotient map $q$ from $\TCP$ to $\CP$ is
covariant for $\tsigma$ and $\csigma$.  Therefore any KMS state $\psi$
on $\CP$ yields the KMS state $\psi\compos q$ on $\TCP$.  Conversely,
any KMS state on $\TCP$ which vanish on $\Ker(q)$ gives a KMS state on
$\CP$ by passage to the quotient.

  Observe that the element $\cS \in\CP$ is analytic for the gauge
action and that
  $$
  \csigma_z(\cS ) = h^{iz}\cS  \for z\in \C.
  $$
  We then have
for every $n\in\N$ that
  $$
  \csigma_{z}(\cS^n) =
  \underbrace{\big(h^{iz}\cS \big) \cdots \big(h^{iz}\cS \big)}_{n\rm\;times} =
  h^{iz} \a\big(h^{iz}\big) \cdots \a^{n-1}\big(h^{iz}\big) \cS^n =
  % \big(h^{iz}\big)^{[n]}
  h^{iz[n]}
  \cS^n,
  $$
  where by   $h^{iz[n]}$ we of course mean $\big(h^{iz}\big)^{[n]}$.
  Since $\csigma_{t}(\cS^*) = \cS^*h^{-it}$, for $t\in\Real$, we have that
$\csigma_{z}(\cS^*) = \cS^*h^{-iz}$, for $z\in\C$, so that for $m\in\N$
  $$
  \csigma_{z}(\cS^m) =
  \underbrace{\big(\cS^*h^{-iz}\big) \cdots \big(\cS^*h^{-iz}\big)}_{m\rm\;times} =
  \cS^{*m} \a^{m-1}\big(h^{iz}\big) \cdots \a\big(h^{iz}\big) h^{iz}  =
  \cS^{*m} h^{-iz[m]}.
  $$
  It is therefore clear that any element of the form
$a\cS^n\cS^{*m}b$, with $a,b\in A$, is analytic and that
  $$
  \csigma_{z}(a\cS^n\cS^{*m}b) =
  ah^{iz[n]} \cS^n \cS^{*m}  h^{-iz[m]} b \for z\in \C.
  $$
  Obviously the same holds for $\tsigma$ and $\tS$.

Our next goal will be the characterization of the states $\phi$ of $A$
such that the composition $\phi \compos \cG$ is a KMS state on $\CP$
(and hence $\phi \compos \cG \compos q$ is a KMS state on $\TCP$).

  \state Proposition
  \label \KMSPhi
  Let $\phi$ be a state on $A$ and let $\b>0$ be a real number.  Then
the state $\psi$ on $\CP$ given by $\psi=\phi \compos \cG$ is a {\kms}
state if and only if $\phi$ is a trace and
  $$
  \phi\big(a\big) = \phi\big(\Tr(\Lambda a)\big)
  $$
  for all $a\in A$, where
  $\Lambda = h^{-\b}\ind(E)$.

  \proof
  Suppose that $\phi$ is a trace satisfying the condition in the
statement.  In view of \lcite{\LinearSpan}, in order to prove $\psi$
to be a {\kms} state it is enough to show that
  $$
  \psi(x\csigma_{i\b}(y)) = \psi(yx)
  \eqno{(\seqnumbering)}
  \label \IWasStar
  $$
  whenever $x$ and $y$ have the form
$x=a\cS^n\cS^{*m}b$ and $y=c\cS^j\cS^{*k}d$, where $n,m,j,k\in\N$ and
$a,b,c,d\in A$.  We have
  $$
  x\csigma_{i\b}(y) =
  a\cS^n\cS^{*m}b\ \csigma_{i\b}(c\cS^j\cS^{*k}d) =
  a\cS^n\cS^{*m}bch^{-\b[j]}\cS^j\cS^{*k}h^{\b[k]}d \&=
  a\a^n\Tr^m(bch^{-\b[j]})\cS^{n-m+j}\cS^{*k}h^{\b[k]}d,
  $$
  by \lcite{\FormulasForLinearSpan} as long as we assume that $m\leq j$.
  We therefore see that $\cG(x\csigma_{i\b}(y)) = 0$ if $n-m+j-k\neq0$.
Under the latter condition it is also easy to see that $\cG(yx) = 0$, in
which case \lcite{\IWasStar} is verified.  We thus assume that $n-m+j-k=0$.
  Setting $p= j-m$ we have that $j=m+p$ and $k=n+p$.  So
  $$
  \psi(x\csigma_{i\b}(y)) =
  \phi\compos \cG\big(a\a^n\Tr^m(bch^{-\b[m+p]})\cS^{n+p}\cS^{*(n+p)}h^{\b[n+p]}d\big) \&=
  \phi\big(a\a^n\Tr^m(bch^{-\b[m+p]})\iota_{n+p}\inv h^{\b[n+p]}d\big) =
  \phi\big(\a^n\Tr^m(u)v\big),
  $$
  where
  $
  u:= bch^{-\b[m+p]}
  $ and $
  v:= \iota_{n+p}\inv h^{\b[n+p]}da.
  $
  On the other hand, 
  $$
  \psi(yx) =
  \psi(c\cS^{m+p}\cS^{*(n+p)}d a\cS^n\cS^{*m}b) \={\FormulasForLinearSpan}
  \phi\compos \cG\big(c\cS^{m+p}\cS^{*(m+p)} \a^m\Tr^n(d a)b \big) =
  \phi \big(c\iota_{m+p}\inv \a^m\Tr^n(d a)b \big) \&=
  \phi \big(uh^{\b[m+p]}\iota_{m+p}\inv \a^m\Tr^n(
    h^{-\b[n+p]} \iota_{n+p} v ) \big) =
  \phi \big(u\Lambda^{-[m+p]} \a^m\Tr^n( \Lambda^{[n+p]}v ) \big) 
  \={\FormuLoca.i}
  %   $$
  %   $$
  %   \=\FormuLoca
  %   \phi \big(u\Lambda^{-[m+p]} \a^m(\Lambda^{[p]})
  %       \a^m\Tr^n( \Lambda^{[n]} v ) \big) =
  \phi \big(u\Lambda^{-[m]} \a^m\Tr^n( \Lambda^{[n]} v ) \big).
  $$
  We thus see that \lcite{\IWasStar} holds under the present
hypotheses that $m\leq j$, if and only if
  $$
  \phi\big(\a^n\Tr^m(u)v\big) =
  \phi \big(u\Lambda^{-[m]} \a^m\Tr^n( \Lambda^{[n]}v ) \big)
  \for u,v\in A \for n,m\in\N.
  \eqno{(\dag)}
  $$
  Consider the linear maps
  $
  \tilde \a, \tilde\Tr : A \to A
  $
  given respectively by
  $\tilde\a(a) = \Lambda\inv \a(a)$,
  and
  $\tilde\Tr(a) = \Tr(\Lambda a)$,
  for all $a\in A$.  It is then easy to see that
  $\tilde\a^m(a) = \Lambda^{-[m]}\a^m(a)$,
  and that
  $\tilde\Tr^n(a) = \Tr^n(\Lambda^{[n]} a)$.  The equation in $(\dag)$
is then expressed as
  $$
  \phi\big(\a^n\Tr^m(u)v\big) =
  \phi \big(u \tilde \a^m\tilde \Tr^n( v ) \big).
  \eqno{(\ddag)}
  $$

  Observe that by hypotheses we have for all $a,b\in A$ that
  $$
  \phi(\a(a)b) =
  \phi(\Tr(\Lambda\a(a)b)) =
  \phi(a\Tr(\Lambda b)) =
  \phi(a\tilde\Tr(b)),
  $$
  and
  $$
  \phi(\Tr(a)b) =
  \phi\big(\Tr(a\a(b))\big) = 
  \phi\big(\Tr(\Lambda \Lambda\inv a\a(b))\big) = 
  \phi\big(\Lambda\inv a\a(b)\big) = 
  \phi\big(a\tilde\a(b)\big).
  $$
  This may be interpreted as saying that with respect to the
inner--product $\<a,b\>=\phi(a^*b)$ one has that the adjoint of $\a$
is $\tilde \Tr$ and the adjoint of $\Tr$ is $\tilde \a$.
  It is now evident that $(\ddag)$ holds, hence completing the proof
of \lcite{\IWasStar}
in the case
that $m\leq j$.

  When $m\geq j$ it is also true that both sides of \lcite{\IWasStar}
vanish unless $n-m+j-k=0$.  In this case, letting $p = m-j$, we have
that $m=j+p$ and $n=k+p$, so that
  $$
  \psi\big(x\csigma_{i\b}(y)\big) =
  \psi\big(a\cS^{k+p}\cS^{*(j+p)}bch^{-\b[j]}\cS^j\cS^{*k}h^{\b[k]}d\big)
  \={\FormulasForLinearSpan}
  \phi\compos \cG\big(a\cS^{k+p}\cS^{*(k+p)}  \a^k\Tr^j(bch^{-\b[j]})h^{\b[k]}d\big) \&=
  \phi\big(a\iota_{k+p}\inv \a^k\Tr^j(bch^{-\b[j]})h^{\b[k]}d\big) =
  \phi\big(\a^k\Tr^j(u)v\big)
  $$
  where
  $
  u =: bch^{-\b[j]}
  $
  and
  $
  v =: h^{\b[k]}d a\iota_{k+p}\inv.
  $
  On the other hand
  $$
  \psi(yx) =
  \psi\big( c\cS^j\cS^{*k}d a\cS^{k+p}\cS^{*(j+p)}b \big)
  \={\FormulasForLinearSpan}
  \phi\compos \cG\big(  c\a^j\Tr^k(d a)\cS^{j+p}\cS^{*(j+p)}b  \big) =
  \phi\big( c\a^j\Tr^k(d a)\iota_{j+p}\inv b \big) \&=
  \phi\big(\a^j\Tr^k(h^{-\b[k]}v\iota_{k+p})\iota_{j+p}\inv uh^{\b[j]}\big) 
  \={\FormuLoca.ii}
  \phi\big(\a^j\Tr^k(h^{-\b[k]}v\iota_{k})\iota_{j}\inv uh^{\b[j]} \big) =
  \phi\big(\acute\a^j\acute\Tr^k(v)u\big),
  $$
  where $\acute\a$ and $\acute\Tr$ are defined respectively by
  $\acute\a(a) = h^\b\a(a)\; \ind(E)\inv$
  and
  $
  \acute\Tr(a) = \Tr(h^{-\b}a\; \ind(E)).
  $
  However, since both $\ind(E)$ and $h$ belong to the center of $A$ we have that
  $\acute\a=\tilde\a$ and   $\acute\Tr=\tilde\Tr$, so that under the
hypotheses that $m\geq j$ we see that  \lcite{\IWasStar} is equivalent to
  $$
  \phi\big(\a^k\Tr^j(u)v\big) =
  \phi\big(u\tilde\a^j\tilde\Tr^k(v)\big),
  $$
  which follows as above.

  Conversely, supposing that $\psi$ is a {\kms} state on $\CP$ we have that
\lcite{\IWasStar} holds for all analytic elements $x$ and $y$.  Given $a,b\in A$ plug
$x=a$ and $y=b$ in \lcite{\IWasStar} to conclude that
  $\phi(ab)=\phi(ba)$ so that $\phi$ must be a trace on $A$.
  On the other hand, plugging $x=\cS^*$ and $y= a\;\ind(E)\cS $ in
\lcite{\IWasStar} gives
  $$
  \phi\big(\Tr(a\;\ind(E)h^{-\b})\big) =
  \phi\big(a\big),
  $$
  hence completing the proof.
  \proofend

  The KMS states provided by the above result necessarily vanish on
elements of the form $a\cS^n\cS^{*m}b$ with $n\neq m$ since so does
$\cG$.  We shall see next that this is necessarily the case for all
KMS states when $h\geq cI$ for some real number $c>1$ (as opposed to
$c>0$ which we have been assuming so far).  We will in fact prove a
slightly stronger result by considering KMS states on $\TCP$, which
include the KMS states on $\CP$ as already mentioned.

  \state Proposition
  \label \VanishLemma
  Suppose that $h\geq cI$ for some real number $c>1$ and let $\psi$ be
a {\kms} state on $\TCP$, where $\b>0$.  Then for every $a,b\in A$ and
every $n,m\in\N$ with $n\neq m$ one has that
$\psi(a\tS^n\tS^{*m}b)=0$.

  \proof
  Taking adjoints we may assume that $n> m$.   So write $n=m+p$ with
$p>0$.
  We have
  $$
  \psi(a\tS^n\tS^{*m}b)=
  \psi(a\tS^m\tS^p\tS^{*m}b)=
  \psi(\tS^p\tS^{*m}b\,\tsigma_{i\b}(a\tS^m)) =
  \psi(\tS^p\tS^{*m}ba h^{-\b[m]}\tS^m) \&=
  \psi(\tS^p \Tr^m(ba h^{-\b[m]})  )=
  \psi(\a^p \Tr^m(ba h^{-\b[m]}) \tS^p ).
  $$
  So it suffices to prove that
  $
  \psi(a\tS^p) =0$ for all $a\in A$ and $p>0$.
  In order to accomplish this notice that
  $$
  \psi(a\tS^p) =
  \psi(\tS^p \tsigma_{i\b}(a)) =
  \psi(\tS^p a) =
  \psi(a\tsigma_{i\b}(\tS^p)) =
  \psi(a h^{-\b[p]}\tS^p),
  $$
  so that 
  $$
  \psi(ak\tS^p)=0 \for a\in A,
  \eqno{(\dag)}
  $$
  where
  $k=1-h^{-\b[p]}$.  Since $h\geq c$ we have that
  $h^{-\b}\leq c^{-\b}$ and hence 
  $h^{-\b[p]}\leq c^{-\b p}$ by \lcite{\CohomologyMiscelanea.iii}.
This implies that
  $k \geq 1 - c^{-\b p} > 0$ and hence that $k$ is invertible.
The
conclusion then follows upon replacing $a$ with  $ak\inv$  in $(\dag)$.
  \proofend

  Observe that we haven't used that $E$ is of index-finite type in the
above proof.  Also notice that it follows from the above result that
any {\kms} state on $\CP$ must vanish on elements of the form
$a\cS^n\cS^{*m}b$ with $n\neq m$.

  We would now like to address the question of whether all KMS state
on $\CP$ are given by \lcite{\KMSPhi}.  Should there exist more than
one conditional expectation from $\CP$ to $A$ it would probably be
unreasonable to expect this to be true.  In view of
\lcite{\ConditionalExpectOnF} and \lcite{\UniqueCondExpInvariant} one
is led to believe that the question posed above is easier to be dealt
with under the hypothesis that $A$ is commutative.

After having proved the result below for  commutative algebras I
noticed that the commutativity hypothesis was used only very slightly
and could be replaced by the weaker requirement that $E(ab)=E(ba)$ for
all $a,b\in A$.  In the hope that a relevant example might be found
under this circumstances we will restrict ourselves to this weaker
hypothesis whenever possible.

  \state Proposition
  \label \MainTechnical
  Suppose that $h\geq cI$ for some real number $c>1$ and let $\psi$ be
a {\kms} state on $\CP$, where $\b>0$.  Suppose moreover 
that $E(ab)=E(ba)$ for all $a,b\in A$ (e.g.~when $A$ is commutative).
Then $\psi = \psi\compos \cG$.  Therefore $\psi$ is given as in
\lcite{\KMSPhi} for $\phi=\psi|_A$.

  \proof
  We shall prove the equivalent statement that all {\kms} states
$\psi$ on $\TCP$ which vanish on $\Ker(q)$ satisfy $\psi = \psi\compos
\tG$.  

  Let $(u_1,\ldots,u_m)$ be a quasi-basis for $E$ as in the beginning
of section \lcite{\FiniteIndexSection}.  Setting
  $k= \sum_{j=1}^m u_j \tS  \tS^* u_j^*$ observe that for all $b\in A$ one
has
  $$
  kb\tS  =
  \sum_{j=1}^m u_j \tS  \tS^* u_j^*b\tS  =
  \sum_{j=1}^m u_j \tS  \Tr(u_j^*b) =
  \sum_{j=1}^m u_j \a (\Tr(u_j^*b)) \tS =
  \sum_{j=1}^m u_j E(u_j^*b) \tS = b\tS ,
  $$
  showing that the pair $(1,k)$ is a redundancy.  It follows that
$1-k\in\Ker(q)$ and hence for all $a\in A$
  $$
  \psi(a) =
  \psi(ak) =
  \psi\(\sum_{j=1}^m au_j \tS  \tS^* u_j^*\) =
  \sum_{j=1}^m \psi(\tS^* u_j^*\tsigma_{i\b}(au_j \tS )) =
  \sum_{j=1}^m \psi(\tS^* u_j^* au_j h^{-\b}\tS ) \&=
  \sum_{j=1}^m \psi(\Tr(u_j^* au_j h^{-\b}))  =
  \sum_{j=1}^m \psi(\Tr( a h^{-\b} u_j u_j^*)) =
  \psi(\Tr(a h^{-\b}\ind(E)))  =
  \psi(\Tr(a\Lambda)),
  $$
  where, as before,   $\Lambda = h^{-\b}\ind(E)$.  Replacing 
$a$ by $a\Lambda\inv$ above leads to
  $\psi(\Tr(a)) =
  \psi(\Lambda\inv a)
  $.  It is then easy to prove by induction that
  $$
  \psi(\Tr^n(a)) =
  \psi(\Lambda^{-[n]}a),
  $$
  for all $a\in A$ and $n\in\N$.
  Given $n,m\in\N$ and $a,b\in A$ we claim that
  $$
  \psi(a\tS^n\tS^{*m}b) =
  \psi(\tG(a\tS^n\tS^{*m}b)).
  $$
  Observe that the case in which $n\neq m$ follows immediately from
\lcite{\VanishLemma}.  So  we assume that $n=m$. We
then have that
  $$
  \psi(a\tS^n\tS^{*n}b) =
  \psi(\tS^{*n}b\tsigma_{i\b}(a\tS^n)) =
  \psi(\tS^{*n}bah^{-\b[n]}\tS^n) =
  \psi(\Tr^n(bah^{-\b[n]})) \&=
  \psi(\Lambda^{-[n]}bah^{-\b[n]}) \=\star
  \psi(ah^{-\b[n]}\Lambda^{-[n]}b) =
  %  \psi(\ind(E)^{-[n]}ba) =
  \psi(a\;\ind(E)^{-[n]}b) =
  \psi(\tG(a\tS^n\tS^{*n}b)),
  $$
  where we have used  in $(\star)$ the fact that the restriction of a KMS state to
the algebra of fixed points is a trace.  This proves our
claim and the result follows from \lcite{\LinearSpan}.
  \proofend

  Summarizing we have:

  \state Theorem
  \label \MainTheorem
  Let $\a$ be an injective endomorphism of a unital C*-algebra $A$
with $\a(1)=1$.  Let $E$ be a conditional expectation of index-finite
type from $A$ onto the range of $\a$ such that $E(ab)=E(ba)$ for all
$a,b\in A$ (e.g.~when $A$ is commutative).  Let $\Tr=\a\inv\compos E$
be the corresponding transfer operator.  Given a self-adjoint element
$h\in\Zenter(A)$ with $h\geq cI$ for some real number $c>1$, consider
the unique one-parameter automorphism group $\csigma$ of $\CP$ given
for $t\in\Real$ by $\sigma_t(\cS) = h^{it}\cS$ and $\sigma_t(a)=a$ for
all $a\in A$.  Then, for all $\b>0$ the correspondence $$\psi\mapsto
\phi=\psi|_A$$ is a bijection from
  the set of {\kms} states $\psi$ on $\CP$ and
  the set of states $\phi$ on $A$ such that
  %  $$
  %  \phi\Big(\Tr\big(\ind(E) a\big)\Big) =
  %  \phi(h^{\b} a) 
  %  \for a\in A.
  %  $$
  $
  \phi\big(a\big) = \phi\big(\Tr(\Lambda a)\big)
  $
  for all $a\in A$, where
  $\Lambda = h^{-\b}\ind(E)$.  The inverse of the above correspondence
is given by $\phi\mapsto\psi=\phi\compos \cG$, where $\cG$ is the
conditional expectation given in \lcite{\UniqueCondExpInvariant}.

  \section{Ground states}
  In this section we retain the standing assumptions made in
\lcite{\StandingHyp} but we will drop \lcite{\StandingHyp.iv} at a
certain point below.  Our goal is to treat the case of ground states
on $\CP$ for the gauge action $\csigma^h$.  Recall that a state $\psi$
on $\CP$ is a ground state if
  $$
  \sup_{\hbox{\sevenrm Im} z \geq 0} |\psi(x\sigma_z(y))| < \infty
  $$
  for every analytic elements $x,y\in\CP$.  Let $(u_1,\ldots,u_m)$ be
a quasi-basis for $E$ as in the beginning of section
\lcite{\FiniteIndexSection}.  As seen in the proof of
\lcite{\MainTechnical} the pair $(1,k)$ is a redundancy, where
  $k= \sum_{j=1}^m u_j \tS  \tS^* u_j^*$.  Therefore one has that
  $$
  1 = \sum_{j=1}^m u_j \cS  \cS^* u_j^* 
  $$
  in $\CP$.  Assuming that $\psi$ is a ground state on $\CP$ one has
that the following is bounded for  $z$ in the upper half
plane:
  $$
  \sum_{j=1}^m \psi(u_j \cS  \csigma_z(\cS^* u_j^*)) =
  \sum_{j=1}^m \psi(u_j \cS  \cS^* h^{-iz} u_j^*) =
  \psi(h^{-iz} ),
  $$
  say by a constant $K>0$.  With $z=i\b$ 
  we conclude  that
  $
  \psi(h^\b) \leq K
  $
  for all $\b>0$.  Suppose that $h\geq cI$ for some real number $c>1$
as before.  Then $h^\b \geq c^\b$ and
  $$
  K \geq \psi(h^\b) \geq c^\b.
  $$
  Observing that the term in right hand side above converges to
infinity as $\b\to\infty$ we arrive at a contradiction thus proving:
 
  \state Proposition
  Suppose that $E$ is of index-finite type and that $h\geq cI$ for
some real number $c>1$.  Then there are no ground states on $\CP$.

In the remainder of this section we will discuss the ground states on
$\TCP$.  Our results in this direction will no longer depend on the
fact that $E$ is of index-finite type.

  \state Proposition
  \label \GroundLemma
  Suppose that $h\geq cI$ for some real number $c>1$.  Then a state
$\psi$ on $\TCP$ is a ground
state if and only if $\psi$ vanishes on any element of the form
  $a \tS^n \tS^{*m}b$
  % , where $a,b\in A$, $n,m\in\N$
  if $(n,m)\neq(0,0)$.

  \proof
  Let $a,b\in A$ and $n,m\in \N$ with $(n,m)\neq(0,0)$ and
let $\psi$ be a ground state on $\TCP$.  By taking adjoints it
suffices to prove the result in the case that $m\neq0$.  Letting
  $x=a \tS^n$ and
  $y=\tS^{*m}b$
  we have that
  $$
  \psi(x\tsigma_z(y)) =
  \psi(a \tS^n\tS^{*m}h^{-iz[m]}b)
  \eqno{(\dag)}
  $$
  is bounded as a function of $z$ on the upper half plane.  For
$z=x+iy$ we have
  $$ 
  \|h^{-iz}\| = 
  \|h^{y-ix}\| = 
  \|h^{y}\|.
  $$
  If $z$ is in the \stress{lower half plane}, that is if $y\leq0$,
then since $h\geq cI$ we have that $h^y \leq c^y < 1$ so that $(\dag)$
is actually bounded everywhere.  By Liouville's Theorem $(\dag)$ is
constant and that constant must be zero since zero is the limit of
$(\dag)$ as $z$ tends to infinity over the negative imaginary axis.
Plugging $z=0$ in $(\dag)$ gives the desired conclusion.
  We leave the proof of the converse statement to the reader.
  \proofend

We now need some insight on the structure of the fixed-point algebra
for the scalar gauge action $\tgamma$ on $\TCP$.

  \state Proposition
  Let $\tU$ be the subalgebra of $\TCP$ consisting of the fixed-points
for $\tgamma$.  Then there exists a *-homomorphism
  $\pi: \tU \to A$
  such that $\pi(a)=a$, for all $a\in A$, and $\pi(\tS^n \tS^{*n})=0$,
for all $n>0$.

  \proof
  Consider the representation $\rho: \TCP\to \Lin(M_\infty)$ described
in the proof of \scite{\Endo}{Theorem 3.4}.
  It is easy to see that $\rho$ maps $\tU$ into the set of diagonal
operators with respect to the decomposition 
  $
  M_\infty = \bigoplus_{n=0}^\infty M_{\Tr^n}.
  $
  Therefore, letting $e$ be the projection onto $M_{\Tr^0}$, we have
that the map
  $$
  \pi:x \in \tU  \longmapsto e\rho(x)e \in\Lin (M_{\Tr^0})
  $$
  is a *-homomorphism.
  It is evident that  $\pi$ maps each $a\in A$ to the same $a$ in
the canonical copy of $A$ within $\Lin(M_0)$ while 
$\pi(\tS^n\tS^{*n}) = 0$ for all $n>0$.
  \proofend
  
  \state Proposition
  Suppose that $h\geq cI$ for some real number $c>1$.  Then 
(regardless of E being of index-finite type or not) the ground
states on $\TCP$ are precisely the states of the form $\phi\compos\pi
\compos \tP$
where 
$\tP$ is the
conditional expectation onto $\tU$ given by
  $$
  \tP(x)= \int_{S^1} \tgamma_z (x)\, dz \for x\in\TCP.
  $$
  and $\phi$ is any state whatsoever on $A$.

  \proof
  Let $\psi$ be a ground state on $A$.  Then as a special case of
\lcite{\GroundLemma} we see that $\psi$ vanishes on $a \tS^n
\tS^{*m}b$ whenever $n\neq m$.  By checking first on the generators of
$\TCP$ provided by \lcite{\LinearSpan} it is easy to see that 
  $\psi = \psi \compos \tP$.  Letting $\chi$ denote the restriction of
$\psi$ to $\tU$ we then evidently have that
  $\psi = \chi \compos \tP$.

  Let now $\phi$ be the restriction of
$\chi$ (and hence also of $\psi$) to $A$.  Then one  may prove that
  $\chi = \phi\compos\pi$ by checking on the generators of $\tU$
given by \lcite{\FixedPointAlgebras}.
  So $\psi = \chi \compos \tP = \phi\compos\pi \compos \tP$ as
desired.

  Conversely, given any state $\phi$ on $A$ it is easy to see that
$\psi = \phi\compos\pi\compos \tP$ is a ground state by
\lcite{\GroundLemma}.
  \proofend
  
  \section{The commutative case}
  Let us now discuss the case of a commutative $A$.  Rather than
employ Gelfand's Theorem and view $A$ as the algebra of continuous
functions on its spectrum we will let $A$ be any closed unital
*-subalgebra of the C*-algebra $B(X)$ of all bounded functions on a
set $X$ (with the sup norm).  Examples are:
  \izitem
  \zitem if $X$ is a measure space take $A$ to be the set of all
bounded measurable functions on $X$,
  \zitem if $X$ is a topological space choose a subset
$\{x_1,x_2,\ldots\}\subseteq X$ and let $A$ be the set of all bounded
functions which are continuous at all points of $X$ except, perhaps,
at the points of the set above.

  \bigskip\noindent
  Let us also fix a surjective mapping
  $$
  \theta : X \to X
  $$
  such that $f\compos \theta \in A$ for all $f\in A$.  Clearly one gets
a unital *-monomorphism
  $
  \a: A \to A
  $
  by letting 
  $$
  \a(f) = f\compos \theta
  \for f\in A.
  $$
  Assume that there exists a finite subset
  $\{v_1, \ldots, v_m\}\subseteq A$ such that for all $i=1,\ldots,m$:
  \izitem
  \zitem $\theta$ is injective when restricted to the set $\{x\in X:
v_i(x)\neq0\}$,
  \zitem $v_i\geq0$,
  \zitem $\sum_{i=1}^m v_i = 1$.

  \bigskip\noindent
  For each $x\in X$ define
  $$
  \numb(x) = \#\big\{t\in X: \theta(t)=x\big\}
  $$
  and observe that the existence of the $v_i$'s above implies that
$\numb(x)\leq m$.  For $f\in A$ consider the function $\T(f)$ on $X$
given by
  $$
  \T(f) \calcat x = 
  \sum_{\buildrel {\scriptstyle t\in X} \over {\theta(t)=x}} f(t).
  $$
  If we assume that $\T(f)\in A$ for all $f\in A$ and moreover that 
$\numb$, seen as a
bounded function on $X$,  belongs to $A$ then the operator
  $\Tr:A\to A$ given by $\Tr(f)= \numb\inv \T(f)$ is a transfer
operator.  In addition the composition $E=\a\compos \Tr$ is a
conditional expectation from $A$ to the range of $\a$, which may be
expressed as 
  $$
  E(f) \calcat x = 
  {1 \over \mu(x)}
  \sum_{\buildrel {\scriptstyle t\in X} \over {\theta(t)=\theta(x)}} f(t)
  $$
  where $\mu=\numb\compos \theta$.
  Setting $u_i=(\mu v_i)^{1/2}$ observe that for all $f\in A$ and
$x\in X$ one has that
  $$
  \sum_{i=1}^m u_i E(u_if)\calcat x =
  \sum_{i=1}^m u_i(x) {1 \over \mu(x)} 
    \sum_{\buildrel {\scriptstyle t\in X} \over {\theta(t)=\theta(x)}}
u_i(t)f(t) =
  \sum_{i=1}^m u_i(x) {1 \over \mu(x)} u_i(x)f(x)=
  \sum_{i=1}^m v_i(x) f(x) = f(x).
  $$
  Therefore $\{u_1,\ldots,u_m\}$ is a quasi-basis for $E$, which says
that $E$ is of index-finite type, and 
  $$
  \ind(E) =
  \sum_{i=1}^m u_i^2 =
  \sum_{i=1}^m \mu v_i = \mu.
  $$
  Fix a positive element $h\in A$ with $h\geq cI$ for some real number
$c>1$ and consider the gauge action $\csigma^h$ on $A$.  By
\lcite{\MainTheorem} we have that the {\kms} states on $\CP$ for the gauge
action $\csigma^h$ correspond to the states $\phi$ on $A$ such that
  $$
  \phi(f) = \phi(\Tr(h^{-\b}\ind(E)f))
  \eqno {(\dag)}
  $$
  for all $f\in A$.  In the present context we have that
  $$
  \Tr(h^{-\b}\ind(E)f)\calcat x =
  {1\over \numb(x)} \sum_{\buildrel {\scriptstyle t\in X} \over
    {\theta(t)=x}} h(t)^{-\b} \mu(t)f(t) =
  \sum_{\buildrel {\scriptstyle t\in X} \over
    {\theta(t)=x}} h(t)^{-\b} f(t).
  $$
  The operator $f\mapsto \Tr(h^{-\b}\ind(E)f)$ therefore coincides
with the operator $L_{h^{-\b}}$ introduced by Ruelle in \cite{\RueleOne},
\cite{\RueleTwo}.

  \endgroup

  \end